\documentclass[12pt]{article}

\author{Stefan Ankirchner\\Institut f\"ur Mathematik\\
Humboldt-Universit\"at zu Berlin\\ Unter den Linden 6\\ 10099
Berlin\\ Germany \and Peter Imkeller\\ Institut f\"ur Mathematik\\
Humboldt-Universit\"at zu Berlin\\ Unter den Linden 6\\ 10099
Berlin\\ Germany \and Alexandre Popier\\Laboratoire manceau de math\'ematiques\\ Universit\'e du Maine\\72085 Le Mans C\'edex 9\\ France }

\title{On measure solutions of backward stochastic differential equations}

\usepackage{amssymb}

\begin{document}

\def\i#1#2#3{\int_{0}^{#1}\, #2\, d{#3}} \def\ii#1#2{\int_{0}^{1} #1\, d{#2}}
\def\s{\sigma}
\def\no#1#2#3{\parallel #1
\parallel_{#2}^{#3}}
\def\a#1{\mid#1\mid}
\def\ind#1{1_{\{ #1\}}}
\def\bn#1{#1\in {\bf N}}
\def\br#1{#1\in {\bf R}}
\def\ui{[0,1]}
\def\uii#1{#1\in \ui}
\def\pf{{\bf Proof: }}
\def\re{{\bf Remark: }}
\def\ep{\epsilon}
\def\al{\alpha}
\def\om{\omega}
\def\Om{\Omega}
\def\iom{\om\in \Om}

\def\f{{\cal F}}
\def\r{{\mathbb{R}}}
\def\n{{\mathbb N}}
\def\z{{\bf Z}}
\def\b{{\bf B}}
\def\d{{\bf D}}
\def\m{{\bf M}}
\def\i{{\cal I}}
\def\l{\lambda}
\def\L{\Lambda}
\def\g{{\cal G}}
\def\ft{f_{1-t}}
\def\un{\infty}
\def\e#1{\exp(#1)}
\def\eh{\fr{1}{2}}
\def\edp{\fr{1}{p}}
\def\ep{\epsilon}
\def\nh{\fr{n}{2}}
\def\p{\phi}
\def\P{\Phi}
\def\px{\p_t(x)}
\def\py{\p_t(y)}
\def\psxs{\ps_t(x,s)}
\def\pih{\fr{\pi}{2}}
\def\lxs{l_t(x,s)}
\def\pxp{\p_t(x\rq)}
\def\psxsp{\ps_t(x\rq, s\rq)}
\def\lxsp{l_t(x\rq,s\rq)}
\def\xp{x^{\rq}}
\def\sp{s^{\rq}}
\def\ps{\psi}
\def\Ps{\Psi}
\def\fr{\frac}
\def\dj{D^j}
\def\ds{D_s}
\def\epf{$\Box$ \\}
\def\st{\sup_{0\leq t\leq 1}}
\def\ss{\sup_{s\in S^{d-1}}}
\def\G{\Gamma}
\def\cc{{\bf C}}
\def\zpi{\fr{1}{\sqrt{2 \pi}}}
\def\ezt{\fr{1}{2t}}
\def\ejm{1\leq j\leq m}
\def\ekr{1\leq k\leq r}
\def\fr#1#2{\frac{#1}{#2}}
\def\tun#1{#1 \to\un}
\def\suim{\sum_{i=1}^m}
\def\sir{\sum_{i=1}^r}
\def\sid{\sum_{i=1}^d}
\def\skd{\sum_{k=1}^d}
\def\skm{\sum_{k=1}^m}
\def\spqd{\sum_{p,q=1}^d}
\def\sqd{\sum_{q=1}^d}
\def\sjm{\sum_{j=1}^m}
\def\spd{\sum_{p=1}^d}
\def\vec#1#2#3{(#1_{#3}^1,...,#1_{#3}^{#2})} \def\su{{\rm supp}}
\def\grad{\mbox{grad}}
\def\sgn{\mbox{sgn}}
\def\der#1{\fr{d}{d#1}}
\def\parder#1#2{\fr{\partial #1}{\partial #2}} \def\vect#1#2{\left[
\begin{array}{cc}
#1 \\
#2
\end{array}
\right]}

\def\matr#1#2#3#4{\left[
\begin{array}{cc}
#1 & #2\\
#3 & #4
\end{array}
\right]}
\def\zixj#1#2#3#4{\left[ \begin{array}{ccc}
#1 & \cdots & #2 \\
\vdots & & \vdots \\
#3 & \cdots & #4
\end{array}
\right]}
\def\nn{\nonumber}
\def\mall#1{\langle\langle D#1,D#1\rangle\rangle} \def\tr{\mbox{tr}}
\def\hpt{\hat{P}_t}
\def\hrt{\hat{R}_t}
\def\la{\langle}
\def\ra{\rangle}

\newtheorem{thm}{Theorem}[section]
\newtheorem{pr}{Proposition}[section]
\newtheorem{co}{Corollary}[section]
\newtheorem{lem}{Lemma}[section]
\newtheorem{defi}{Definition}[section]
\newtheorem{rem}{Remark}[section]

\newcommand{\eps}{\varepsilon}
\newcommand{\R}{\mathbb{R}}
\newcommand{\N}{\mathbb{N}}
\newcommand{\E}{\mathbb{E}}
\newcommand{\prb}{\mathbb{P}}
\newcommand{\qprb}{\mathbb{Q}}
\newcommand{\bp}{{\bf P}}
\newcommand{\expo}{\mathcal{E}}
\newcommand{\be}{\begin{equation}}
\newcommand{\ee}{\end{equation}}
\newcommand{\bdm}{\begin{displaymath}}
\newcommand{\edm}{\end{displaymath}}
\newcommand{\bean}{\begin{eqnarray}}
\newcommand{\eean}{\end{eqnarray}}
\newcommand{\bea}{\begin{eqnarray*}}
\newcommand{\eea}{\end{eqnarray*}}
\newcommand{\zi}[1]{(\ref{#1})}
\newcommand{\wtilde}{\widetilde}

\newcommand{\tri}{\mathcal{F}}
\newcommand{\hilb}{\mathcal{H}}
\newcommand{\ga}{\gamma}
\newcommand{\La}{\Lambda}
\newcommand{\phii}{\varphi}

\newcommand{\cA}{\mathcal{A}}
\newcommand{\cB}{\mathcal{B}}
\newcommand{\cC}{\mathcal{C}}
\newcommand{\cD}{\mathcal{D}}
\newcommand{\cE}{\mathcal{E}}
\newcommand{\cF}{\mathcal{F}}
\newcommand{\cG}{\mathcal{G}}
\newcommand{\cH}{\mathcal{H}}
\newcommand{\cI}{\mathcal{I}}
\newcommand{\cJ}{\mathcal{J}}
\newcommand{\cK}{\mathcal{K}}
\newcommand{\cL}{\mathcal{L}}
\newcommand{\cM}{\mathcal{M}}
\newcommand{\cN}{\mathcal{N}}
\newcommand{\cO}{\mathcal{O}}
\newcommand{\cP}{\mathcal{P}}
\newcommand{\cQ}{\mathcal{Q}}
\newcommand{\cR}{\mathcal{R}}
\newcommand{\cS}{\mathcal{S}}
\newcommand{\cT}{\mathcal{T}}
\newcommand{\cU}{\mathcal{U}}
\newcommand{\cX}{\mathcal{X}}
\newcommand{\cY}{\mathcal{Y}}

\maketitle

\begin{abstract}
We consider backward stochastic differential equations (BSDE) with
nonlinear generators typically of quadratic growth in the control
variable. A measure solution of such a BSDE will be understood as
a probability measure under which the generator is seen as
vanishing, so that the classical solution can be reconstructed by
a combination of the operations of conditioning and using
martingale representations. In case the terminal condition is bounded
and the generator fulfills the usual
continuity and boundedness conditions, we show that measure solutions
with equivalent measures just reinterpret classical ones. In case of
terminal conditions that have only exponentially bounded moments, we
discuss a series of examples which show that in case of non-uniqueness
classical solutions that fail to be measure solutions can coexists with
different measure solutions.

\end{abstract}



\section*{Introduction}

The generally accepted natural framework for the most efficient
formulation of pricing and hedging contingent claims on complete
financial markets, for instance in the classical Merton-Scholes
problem, is given by martingale theory, more precisely by the
elegant notion of martingale measures. Martingale measures
represent a view of the world in which price dynamics do not have
inherent trends. From the perspective of this world, pricing a
claim amounts to taking expectations, while hedging boils down to
pure conditioning and using martingale
representation.\par\smallskip

At first glance, hedging a claim is, however, a problem calling
upon stochastic control: it consists in choosing strategies to
steer the portfolio into a terminal random endowment the portfolio
holder has to ensure. Solving stochastic backward equations (BSDE)
is a technique tailor-made for this purpose. This powerful tool
has been introduced to stochastic control theory by Bismut
\cite{bismut76}. Its mathematical treatment in terms of stochastic
analysis was initiated by Pardoux and Peng \cite{pardoux90}, and
its particular significance for the field of utility maximization
in financial stochastics clarified in El Karoui, Peng and Quenez
\cite{elkarouipengquenez97}. To fix ideas, we restrict our
attention to a Wiener space probabilistic environment. In this
framework, a BSDE with terminal variable $\xi$ at time horizon $T$
and generator $f$ is solved by a pair of processes $(Y,Z)$ on the
interval $[0,T]$ satisfying

\be\label{eq:introduction} Y_t = \xi -
\int_t^T Z_s d W_s + \int_t^T f(s, Y_s, Z_s) ds, \quad
t\in[0,T].\ee
In the case of vanishing generator, the solution
just requires an application of the martingale representation
theorem in the Wiener filtration, and $Z$ will be given as the
stochastic integrand in the representation, to which we will refer
as {\em control process} in the sequel. The classical approach of
existence and uniqueness for BSDE involves a priori inequalities
as a basic ingredient, by which unique solutions are constructed
via fixed point arguments, just as in the case of forward
stochastic differential equations.
\par\medskip

In this paper we are looking for a notion in the context of BSDE that plays the role of the martingale measure
in the context of hedging claims. Our main interest is directed to BSDE of the type (\ref{eq:introduction})
with generators that are non-Lipschitzian, and depend on the control variable $z$ quadratically,
typically $f(s,y,z) = z^2\, b(s,z), s\in[0,T], z\in\r,$ with a bounded function $b$. These generators
were given a through treatment in Kobylanski \cite{kobylanski00}, Briand, Hu \cite{briandhu05}, and Lepeltier,
San Martin \cite{lepeltiersanmartin98}. While \cite{kobylanski00} and \cite{lepeltiersanmartin98} consider
existence and uniqueness questions for bounded terminal variables $\xi$, \cite{briandhu05} goes to the limit
of possible terminal variables by considering $\xi$ for which $\exp(\gamma|\xi|)$ has finite expectation for
some $\gamma > 2 ||b||_{\infty}$. All these papers employ different methods of approach following the classical
pattern of arguments mentioned above. In contrast to this, we shall investigate an alternative notion of solution
of BSDE, the generators of which fulfill similar conditions. In analogy with martingale measures in hedging which
effectively eliminate drifts in price dynamics, we shall look for probability measures under which the generator
of a given BSDE is seen as vanishing. Given such a measure $\qprb$ which we call {\em measure solution} of the
BSDE and supposing that $\qprb \sim \prb$, the processes $Y$ and $Z$ are the results of projection and
representation respectively, i.e. $Y = \E^{\qprb} (\xi|\f_\cdot) = Y_0 + \int_0^\cdot Z_s d \wtilde{W}_s,$
where $\wtilde{W}$ is a Wiener process under $\qprb$. The first main finding of the paper roughly states that provided the
terminal variable $\xi$ is bounded, all classical solutions can be interpreted as measure solutions. More precisely,
we show that if the generator satisfies the usual continuity and quadratic boundedness conditions,
classical solutions $(Y,Z)$ exist if and only if measure solutions with $\qprb \sim \prb$ exist. So
existence Theorems obtained in the papers quoted are recovered in a more elegant and concise way
in terms of measure solutions. We do not touch uniqueness questions in general. Of course,
determining a measure $\qprb$ under which the generator vanishes amounts to doing a Girsanov change of probability
that eliminates it. We therefore have to look at the BSDE in the form
\be\label{eq:introduction2}
Y_t = \xi - \int_t^T Z_s \left[ d W_s - \frac{f(s, Y_s, Z_s)}{Z_s} ds \right],\quad t\in[0,T],\ee
define $\displaystyle g(s,y,z) = \frac{f(s,y,z)}{z}$, and study the measure
$$\qprb = \exp \left( M - \eh \langle M\rangle \right) \cdot \prb$$
for the martingale $M = \int_0^\cdot g(s, Y_s, Z_s) d W_s.$ One of the fundamental problems that took
some effort to solve consists in showing that $\qprb$ is a probability measure. Here one has to dig
essentially deeper than Novikov's or Kazamaki's criteria allow. We successfully employed a criterion
which is based on the explosion properties of the quadratic variation $\langle M\rangle$, which we
learnt from a conversation with M. Yor, and has been latent in the literature for a while, see Liptser,
Shiryaev \cite{liptsershiryaev01}, or the more recent paper by Wong, Heyde \cite{wongheyde04}. This criterion
allows a simple treatment of the problem of existence of measure solutions in the case of bounded
terminal variable, and a still elegant and efficient one in the borderline case of exponentially integrable
terminal variable considered by Briand, Hu \cite{briandhu05}. If $\xi$ is only exponentially bounded,
things turn essentially more complex immediately. Specializing to a very simple generator, we find a wealth of different
situations looking confusing at first sight.
Just to quote three basic scenarios exhibited in a series of examples of different types:
in the first type we obtain one solution which is a measure solution at the same time;
in the second one we find two different solutions both of which are measure solutions;
in the third one we encounter two solutions one of which is a measure solution, while the other one
is not. We even combine these basic examples to develop a scenario in which there exists a continuum of
measure solutions, and another one in which a continuum of non-measure solutions is given.
\par\bigskip

Here is an outline of the presentation of our material. Throughout
we consider BSDE possessing generators with quadratic nonlinearity
in $z$. In a first section we discuss the case of bounded terminal
variable $\xi$, and show that if the generator satisfies continuity
and quadratic boundedness conditions, classical solutions $(Y,Z)$
exist if and only if measure solutions with $\qprb \sim \prb$ exist.
Things become essentially more complex in the second section, where
we pass to exponentially integrable terminal variables. Taking the
simple generator $f(s,z) = \alpha z^2, s\in[0,t], z\in \r,$ with
some $\alpha\in\r$, a wealth of different scenarios arises in which
in case of non-uniqueness in particular solutions can be measure
solutions, while different ones fail to have this property. In the
terminal section we construct measure solutions from first
principles without using strong solutions in our algorithm, for
generators which are Lipschitz continuous with time dependent and
random constants. By iterating the successive applications of
martingale representation and Girsanov change of measure with
respect to drifts obtained from the martingale representation
density of the previous step we obtain a sequence of probability
measures which can be seen to be tight in the weak topology, and
thus have accumulation points which yield measure solutions.

\section{Measure solutions: Definition and first examples}\label{sec1}
In this section we first recall some basic definitions concerning BSDEs. We then introduce and exemplify the notion of a measure solution by looking at a special class of BSDEs.

Throughout let $T$ be a non-negative real, $(\Omega,
\cF, \prb)$ a probability space, and $(W_{t})_{0 \leq t \leq T}$ a one-dimensional
Brownian motion, whose natural filtration, augmented by
$\mathcal{N}$, is denoted by $(\cF_{t})_{0 \leq t \leq T}$, where
$$\mathcal{N} =
\left\{ A \subset \Omega, \ \exists G \in \cF, \ A \subset G \ \mbox{and} \
\prb(G)=0 \right\}.$$

Let $\xi$ be an $\cF_T$-measurable random variable, and let $f: \Omega \times[0,T] \times \R \to \R$ be a measurable function such that for all $z \in \R$ the mapping $f(\cdot, \cdot, z)$ is predictable. A classical solution of the BSDE with {\em terminal condition} $\xi$
and {\em generator} $f$ is defined to be a pair of predictable processes
$(Y,Z)$ such that almost surely we have $\int_0^T Z_s^2 ds < \infty$, $\int_0^T |f(s,Y_s, Z_s)| ds < \infty$, and for all $t\in[0,T]$,
\be\label{BSDE}
 Y_t = \xi - \int_t^T Z_s d W_s + \int_t^T f(s, Z_s) ds.
\ee
The solution processes $(Y,Z)$ are often shown to satisfy some integrability properties and to belong to the following function spaces. For $p \ge 1$ let $\cH^p$ denote the set of all $\R$-valued predictable
processes $\zeta$ such that $E\int_0^1 |\zeta_t|^p dt < \infty$,
and by 
$\cS^\infty$ we denote the
set of all essentially bounded $\R$-valued predictable processes.

If $\xi$ is square integrable and $f$ satisfies a Lipschitz
condition, then it is known that there exists a unique pair $(Y,Z)
\in \cH^2\otimes\cH^2$ solving (\ref{BSDE}). Recall that the
solution process $Y_t$ has a nice representation as a conditional
expectation with respect to a new probability measure if $f$ is a
linear function of the form \be\label{lipbed} f(s,z) = b_s z, \ee
where $b$ is a predictable and bounded process. More precisely, if
$D_t = \exp(\int_0^t b_s dW_s - \frac12 \int_0^t b^2_s ds)$, and
$\mathbb{Q}$ is the probability measure with density $\mathbb{Q} =
D_T \cdot \mathbb{P}$, then \be Y_t = \E^\mathbb{Q}[\xi | \cF_t].
\ee In the following we will discuss whether $Y$ still can be
written as a conditional expectation of $\xi$ if $f$ does not have a
representation as in (\ref{lipbed}) with $b$ bounded, but satisfies
only a quadratic growth condition in $z$. We aim at finding
sufficient conditions guaranteeing that the process $Y_t$ of a
classical solution of a quadratic BSDE has a representation as a
conditional expectation of $\xi$ with respect to a new probability
measure. For this purpose we consider the class of generators
$f:\Omega\times [0,T]\times \r \to \r,$ satisfying for some constant
$c\in\r_+$,

{\flushleft \bf Assumption (H1)}:
\begin{itemize}
\item[(i)] $f(s,z) = f(\cdot,s,z)$ is adapted for any $z\in\r$,
\item[(ii)] $g(s,z)= \frac{f(s,z)}{z}$, $z\in\R$, is continuous in $z$, for all $s\in[0,T]$,
\item[(iii)] $|f(s,z)|\le c(1+z^2)$ for any $s\in[0,T], z\in \R.$
\end{itemize}
Let $\xi$ be an $\cF_T-$ measurable random variable. We introduce for BSDEs with generators satisfying (H1) our concept of measure solutions.
\begin{defi}
A triplet $(Y,Z,\qprb)$ is called {\em measure solution} of the BSDE (\ref{BSDE}), if $\qprb$ is a probability measure on $(\Om, \f)$, $(Y,Z)$ a pair of $(\f_t)$--predictable stochastic processes such that $\int_0^T Z_s^2 ds < \infty$, $\qprb$-a.s. and the following conditions are satisfied:
\begin{eqnarray*}
\wtilde{W} &=& W -  \int_0^\cdot g(s, Z_s) ds\quad\mbox{is a}\quad \qprb-\ \mbox{Brownian motion},\\
\xi &\in& L^1(\Om,\f, \qprb),\\
Y_t &=& \E^\qprb(\xi|\f_t) = \xi - \int_t^T Z_s d \wtilde{W}_s,\quad
t\in[0,T].
\end{eqnarray*}
\end{defi}

It is known from the literature that if the terminal condition $\xi$ is bounded and the generator $f$ satisfies Assumption (H1), then the BSDE (\ref{BSDE}) has a classical solution $(Y,Z)$ (see for example Kobylanski \cite{kobylanski00}). We show that in this case there exists a probability measure $\qprb$, equivalent to $\prb$, such that $(Y,Z,\qprb)$ is a measure solution.
\begin{thm}\label{existence_bounded_liptser}
Assume that $\xi$ is bounded, and that $f$ satisfies Assumption (H1). Then for every classical solution $(Y,Z)$ there exists a probability measure $\qprb$, equivalent to $\prb$, such that $(Y,Z, \qprb)$ is a measure solution of (\ref{BSDE}).
\end{thm}
\pf Let $(Y,Z)$ be a classical solution of (\ref{BSDE}). The very definition entails
$$\int_0^T Z_s^2 ds < \infty, \quad \prb-a.s.$$
Note that due to (ii) and (iii)
$$|g(s,z)|^2 \le C(1+z^2),\quad\mbox s\in[0,T], z\in\r^d,$$
for some $C > 0 $, and hence we have
\be\label{square_integrable} \int_0^T g^2(s, Z_s) ds < \infty, \quad \prb-a.s.\ee
We shall prove that under this condition also a measure solution exists. For this purpose, we define
\be \label{martingaleM}
M =  \int_0^\cdot g(s,Z_s) dW_s.
\ee
It is clear that all we have to establish is that the measure
$$\qprb = V_T\cdot \prb,$$
with
$$V = \exp \left( M - \eh \langle M\rangle \right)$$
leads to a probability measure equivalent to $\prb$. This will be done by investigating possible explosions of the quadratic variation process $\langle M\rangle.$ For $n\in\n$, let
\be \label{explosiontime} \tau_n = T \wedge \inf \{ t\ge 0: \langle M\rangle_t \ge n\}. \ee
Let
$$\mathbb{Q}^n = V_T \cdot \prb |_{\f_{\tau_n}}$$
be the measure change defined locally on $\f_{\tau_n}.$ We know that
$\mathbb{Q}^n$ is a probability measure equivalent to $\prb$, and
the Radon-Nikodym density of $\mathbb{Q}^n$ with respect to $\prb$
on $\f_{\tau_n}$ is given by
$$V_{\tau_n} = \exp \left( M_{\tau_n} - \eh \langle M\rangle_{\tau_n} \right).$$
Moreover, the drifted process
$$\wtilde{W}^n = W - \int_0^{\tau_n\wedge\cdot} g(s, Z_s) ds$$
is a $\mathbb{Q}^n$-- Brownian motion, in particular locally up to
time $\tau_n.$ In order to show that $\qprb$ is a probability
measure, it is sufficient to verify \be\label{explosion}
\mathbb{Q}^n(\tau_n < T) \to 0 \quad (n\to\infty). \ee Namely,
(\ref{explosion}) implies \be \lim_n \E(V_T 1_{\{\tau_n = T\}}) =
\lim_n [\E(V_{\tau_n}) - \E(V_{\tau_n} 1_{\{\tau_n < T\}})] = 1 -
\lim_n \mathbb{Q}^n(\tau_n < T) = 1. \ee On the other hand,
dominated converges yields $E(V_T) = \lim_n \E(V_T 1_{\{\tau_n =
T\}})$, and hence that $\qprb$ is a probability measure. We remark
that the criterion (\ref{explosion}) can be found in
\cite{liptsershiryaev01}, and appears also as Lemma 1.5 in
\cite{Kaz}.

Recall that by the very definition of the measure change,
$$Y_{\tau_n \wedge \cdot} = Y_0 + \int_{0}^{\tau_n \wedge \cdot} Z_s d \wtilde{W}^n_s$$
is a martingale under $\mathbb{Q}^n$, up to time $\tau_n$, which is
bounded by a constant $c_1$, due to the boundedness of $\xi$ (see
Theorem 2.3 in \cite{kobylanski00}).

Hence we obtain for any $n\in\n$, starting with an application of
Chebyshev-Markov's inequality, and, due to (iii), another constant
$c_2$ independent of $n$, such that \bea \mathbb{Q}^n(\tau_n< T) &
\le & \frac{1}{n}\, \E^{\mathbb{Q}^n} \left( \int_0^{\tau_n} g(s,
Z_s)^2 ds \right) = \frac{1}{n}\, \E^n \left( \int_0^{\tau_n} g(s,
Z_s)^2 ds \right) \\\nn &\le& \frac{1}{n} \, c_2 \left( 1+ \E^n
\int_0^{\tau_n} Z_s^2 ds \right) \\\nn &=& \frac{1}{n} \, c_2 \left(
1+  \E^n \left| \int_0^{\tau_n} Z_s d \tilde{W}^n_s \right|^2
\right) = \frac{1}{n} \, c_2 \left( 1+ \E^n \left| Y_{\tau_n} - Y_0
\right|^2 \right) \\\nn &\le& \frac{1}{n}\, c_2 (1 + c_1). \eea Thus
we have shown (\ref{explosion}), and hence that $\qprb$ is a
probability measure. Under $\qprb$, by definition,
$$W^\qprb = W - \int_0^\cdot g(s, Z_s)\, ds$$
is a Brownian motion, and our BSDE may be written as
$$Y_t = \xi - \int_t^T Z_s d W^\qprb_s = \E^\qprb (\xi|\f_t)$$
for $t\in[0,T].$ This shows that $(Y,Z,\qprb)$ is a measure solution. \hfill \epf

It is straightforward to see that every measure solution gives rise to a classical solution. Consequently, under the assumptions of Theorem \ref{existence_bounded_liptser}, measure solutions exist if and if only classical solutions exist. More precisely, we obtain the following.
\begin{co}\label{existence_bounded_liptser_co}
Assume that $\xi$ is bounded, and that $f$ satisfies Assumption (H1). Then $(Y,Z)$ is a classical solution if and only if there exists a probability measure $\qprb$, equivalent to $\prb$, such that $(Y,Z, \qprb)$ is a measure solution of (\ref{BSDE}).
\end{co}

We remark that the previous results can be extended to the case
where $W$ is a $d$-dimensional Brownian motion. Let $f:\Omega\times
[0,T]\times \r^d \to \r$ be a generator for which there exists a
constant $c\in\r_+$ such that \be |f(s,z)| \le c(1+|z|^2),\qquad
s\in[0,T], z\in \R^d, \ee and assume that $g: \Omega \times [0,T]
\times \R^d \to \R^d$ is a function that is continuous in $z$ and
satisfies \be\label{scalarprod} \langle z, g(s,z) \rangle =
f(s,z),\quad \textrm{ for all }z\in \R^d \textrm{ and }s \in [0,T].
\ee If $\xi$ is bounded and $\cF_T$-measurable, then one can show
with similar arguments as used in the preceding proof that, starting
from a classical solution $(Y,Z)$, there exists a probability
measure $\qprb$ such that $W - \int_0^\cdot g(s,Z_s)\, ds$ is a
$\qprb$-Brownian motion, and $Y_t = E^\qprb(\xi|\cF_t)$.

Notice that the relation (\ref{scalarprod}) may be satisfied by more than one continuous $g$, and consequently there may exist more than one measure solution in the multidimensional case. For example, let $d = 2$, $f(s,z) = z_1 z_2$, and observe that $|f(z)| \le \frac12 |z|^2$. For any $a\in(0,\infty)$ let $g_a(z) = (a z_1, \frac1a z_2)$. Then, we have $\langle z, g_a(s,z) \rangle = f(s,z)$, and thus there exist more than one measure solution for a BSDE with generator $f$ and a bounded terminal condition $\xi$.

In the following sections we shall discuss quadratic BSDEs with terminal conditions that are not bounded.
As is known from literature, see for example
Briand, Hu \cite{briandhu05}, \cite{briandhu07}, this case is by
far more complex. For example, it is here that even if the
generators are smooth, solutions stop to be unique. We shall
exhibit examples below which complement the result shown in
Briand, Hu \cite{briandhu07}, according to which uniqueness is
granted in case the generator of the BSDE possesses additional
convexity properties, and the terminal variable possesses
exponential moments of all orders. This fact underlines that also
variations in the generator affect questions of existence and
uniqueness of solutions a lot. For this reason, and also to keep
better oriented on a windy track with many bifurcations, in the next section we shall
choose a simpler generator, and assume that our generator is given by
$$f(s,z) =  \alpha z^2.$$

\section{Measure and non-measure solutions of quadratic BSDEs with unbounded terminal condition}


In this section we will study in more detail the BSDEs with generator of the form
$$f(z) =  \alpha z^2.$$
We shall further assume without loss of generality that $\alpha > 0.$ This can always be obtained in our BSDE by changing the signs
of $\xi$, and the solution pair $(Y,Z).$

Nonetheless, it turns out that positive and negative terminal variables need a separate treatment. We will first show (see Subsection \ref{subsection:lowerexponential}) the existence of measure solutions for terminal conditions $\xi$ bounded from below. Note that by a linear shift of $Y$ we may assume that $\xi\ge 0$. We shall further work under exponential integrability assumptions in the spirit of Briand, Hu \cite{briandhu05}. According to this paper, exponential
integrability of the terminal variable of the form
\be\label{eq:exp_boundedness} \E(\exp(\gamma |\xi|))< \infty\ee
for some $\gamma > 2\alpha$ is sufficient for the existence of a solution. Let us first exhibit an example to show that one cannot go essentially
beyond this condition without losing solvability.

{\bf Example:}\\
Let $T=1,$ and let $\alpha = \eh.$ Let us first consider
$$\xi = \frac{W_1^2}{2}.$$
It is immediately clear from the fact that $W_1$ possesses the standard normal density, that $\E \exp( 2 \alpha |\xi| ) =
\infty,$ hence of course also for $\gamma > 2 \alpha$ (\ref{eq:exp_boundedness}) is not satisfied. To find a solution $(Y,Z)$ of (\ref{BSDE}) on any interval $[t,1]$ with $t>0$ define
$$Z_s = \frac{W_s}{s},\quad s>0,$$
and set for completeness $Z_0 = 0.$ Let $t>0$ and use the product formula for It\^o integrals to deduce
\bean\label{e:intbyparts}
\int_t^1 Z_s dW_s &=& \eh \frac{W_s^2}{s}|_t^1 + \eh \int_t^1
\frac{W_s^2}{s^2} ds\\\nn &=& \xi - \eh \frac{W_t^2}{t} + \eh
\int_t^1 Z_s^2 ds. \eean
This means that, if we set for convenience again $Y_0 = 0$, the pair of processes $(Y_s, Z_s) = (\eh \frac{W_s^2}{s}, \frac{W_s}{s}), s\in\ui,$ solves the BSDE (\ref{BSDE}) on $[t,1]$ for any $t>0.$ Of course, the definition of $Y_0$ is totally inconsistent with the BSDE. Worse than that,
$Z$ is not square integrable on $[0,1]$, as is well known from the path behavior of Brownian motion. Hence $(Y,Z)$ is not a solution of (\ref{BSDE}). To put it more strictly, there is no classical solution of (\ref{BSDE}) on $[0,1]$, since, due to local Lipschitz conditions, any such solution would have to coincide with $(Y,Z)$ on any interval $[t,1]$ with $t>0.$
\par\medskip

According to Jeulin, Yor \cite{jeulinyor90}, transformations of this type are related to a phenomenon they call {\em appauvrissement de filtrations}. In fact, if $\eh$ is replaced with a parameter $\lambda$, they show that the natural filtration of the transformed process gets poorer than the one of the Wiener process, iff $\lambda > \eh.$ Hence in the case we are interested in the Wiener filtration is preserved.
\par\smallskip

Let us now reduce the factor of $W_1^2$ in the definition of $\xi$ a bit, to show that solutions exist in this setting. For $k\in\n,$ let
$$\xi_k = \frac{W_1^2}{2(1+ 1/k)},$$
and consider the BSDE (\ref{BSDE}) with the generator $f$ chosen above, and terminal condition $\xi_k.$ In this setting, we clearly have
$$\E \exp(\gamma \xi_k ) < \infty\quad\mbox{for}\quad 2\alpha \le \gamma < 2\alpha (1+1/k).$$
This shows that the condition of Briand, Hu \cite{briandhu05} is satisfied. It is not hard to construct the solutions of the corresponding BSDEs explicitly, in the same way as above. In fact, for $k\in\n$ we may define $f_k(t) = \frac{1}{k}+t, t\in\ui,$ and set
$$Z^k_t = \frac{W_t}{f_k(t)}, \quad t\in\ui.$$
We may then repeat the product formula for It\^o integrals argument used above to obtain for $t\ge 0$
\bean\label{e:intbyparts_k} \int_t^1 Z^k_s dW_s &=& \eh
\frac{W_s^2}{f_k(s)}|_t^1 + \eh \int_t^1 \frac{W_s^2
f_k'(s)}{f_k(s)^2} ds\\\nn &=& \eh \frac{W_1^2}{f_k(1)} - \eh
\frac{W_t^2}{f_k(t)} + \eh \int_t^1 (Z^k_s)^2 ds. \eean
Hence we set
$$Y_t^k = \eh \frac{W_t^2}{f_k(t)},\quad t\in\ui,$$
to identify the pair of processes $(Y^k, Z^k)$ as a solution of the BSDE
\be\label{BSDE_k} Y^k_t = \xi_k - \int_t^1 Z^k_s d W_s + \eh \int_t^1 (Z^k_s)^2 ds,\quad t\in\ui. \ee
We do not know at this moment whether (\ref{BSDE}) possesses more solutions.\hfill \epf

\subsection{Exponentially integrable lower bounded terminal variable } \label{subsection:lowerexponential}

Under the exponential integrability assumption $\E(\exp(2\alpha \xi))<\infty$, we will now
derive measure solutions from given classical solutions. Leaving the
difficult question of uniqueness apart for a moment, we remark
that with our simple generator, we obtain an explicit solution
given by the formula \be\label{solution_explicit} Y_t =
\frac{1}{2\alpha}\ln M_t - \frac{1}{2\alpha} \ln M_0,\quad Z_t =
\frac{1}{2\alpha} \frac{H_t}{M_t}, \ee where
$$M_t = \E(\exp(2\alpha \xi)|\f_t) = M_0 + \int_0^t H_s d W_s,\quad t\in[0,T].$$
In the sequel, we shall work with this explicit solution. In the
following lemma, we prove integrability properties for the square
norm of $Z$ which will be crucial for stating the martingale
property of $M$ and other related processes later.

\begin{lem}\label{lem:moments}
For any $p\ge 1$ we have
$$\E \left( \left[ \int_0^T Z_s^2 ds\right]^p \right) < \infty.$$
In particular, $\int_0^\cdot Z_s d W_s$ is a uniformly integrable martingale.
\end{lem}

\pf Let $t\in[0,T].$ By It\^o's formula, applied to $N$
$$ \frac{1}{2\alpha}[\ln M_t - \ln M_0] = \frac{1}{2\alpha}\left[ \int_0^t \frac{H_s}{M_s} d W_s - \eh \int_0^t \left( \frac{H_s}{M_s} \right)^2 ds \right] = \int_0^t Z_s d W_s - \alpha \int_0^t Z_s^2 ds.$$
Hence
\be\label{ito}\alpha \int_0^t Z_s^2 ds = -\frac{1}{2\alpha}[\ln M_t - \ln M_0] + \int_0^t Z_s d W_s.\ee
By concavity of the $\ln$ and Jensen's inequality
$$\ln M_t = \ln \E(\exp(2\alpha \xi)|\f_t) \ge \E(2\alpha \xi|\f_t).$$
Using this in (\ref{ito}), we obtain
$$\alpha \int_0^t Z_s^2 ds \le - \E(\xi|\f_t) + \frac{1}{2\alpha} \ln M_0 + \int_0^t Z_s d W_s.$$
Taking $p-$norms in this inequality and using the inequality of Burkholder, Davis and Gundy for the stochastic integral, we obtain with universal constants $c_1, c_2, c_3$
\bea
\E \left( \left[ \int_0^t Z_s^2 ds \right]^p \right) &\le& c_1 \left[ \E \left( |\E(\xi|\f_t)|^p \right) + |\ln M_0|^p + \E \left( \left[ \int_0^t Z_s^2 ds \right]^{\frac{p}{2}} \right) \right]\\
&\le& c_2 \left[ \E(|\xi|^p) + |\ln M_0|^p + \E \left( \left[ \int_0^t Z_s^2 ds\right]^{\frac{p}{2}} \right) \right].
\eea
By a standard argument this entails
$$\E\left( \left[ \int_0^t Z_s^2 ds \right]^p \right) \le c_3[\E(|\xi|^p) + |\ln M_0|^p + 1],$$
and finishes the proof. \hfill \epf

We shall now prove that $(Y,Z)$ gives rise to a measure solution.

\begin{thm}\label{existence_bounded_below}
Assume that $(Y,Z)$ are defined as in (\ref{solution_explicit}). Then there exists a probability measure $\qprb$, equivalent to $\prb$, such that $(Y,Z,\qprb)$ is a measure solution of (\ref{BSDE}).
\end{thm}

\pf Let
$$S = \int_0^\cdot Z_s \, dW_s.$$
Due to Lemma \ref{lem:moments}, we know that $S$ is a uniformly integrable martingale. We may write
\bean\label{Ito_2} \alpha S - \eh \alpha^2 \langle S\rangle &=& \alpha \left[ \int_0^\cdot Z_s d W_s - \alpha \int_0^\cdot Z_s^2 ds \right] +
\int_0^\cdot (\alpha^2 Z_s^2 - \eh \alpha^2 Z_s^2)\, ds\\\nn &=& \alpha(Y - Y_0) + \eh \alpha^2 \int_0^\cdot  Z_s^2 ds. \eean
Now define stopping times $\tau_n=T \wedge \inf\{t\ge 0: \langle S\rangle_t \ge n \}$. For any $n\in\n$ we have
$$\E \exp \left(\alpha S_{\tau_n} - \eh \alpha^2 \langle S\rangle_{\tau_n} \right) =  1,$$
and consequently Fatou's lemma implies
\be\label{Fatou}
\E \exp \left( \alpha [\xi - Y_0] + \eh \alpha^2 \int_0^T Z_s^2 ds\right)\leq \liminf_{n\to\infty} \E \exp \left( \alpha S_{\tau_n} - \eh\alpha^2 \langle S\rangle_{\tau_n} \right) = 1.
\ee
Using this and the positivity of the terminal variable $\xi$, we can now obtain the exponential integrability property
\be\label{exponential_3} \E\exp \left[ \eh \alpha (\xi - Y_0) + \eh \alpha^2 \int_0^T Z_s^2\, ds \right] < \infty. \ee
We shall now use (\ref{Ito_2}) together with (\ref{Fatou}) to prove the exponential integrability of $\eh \alpha S_T$. In fact, we have
$$ \eh \alpha S_T =  \eh\alpha (\xi - Y_0) + \eh \alpha^2 \int_0^T Z_s^2\, ds.$$
Hence we obtain
\be\label{kazamaki} \E \exp \left( \eh \alpha S_T \right) < \infty, \ee
and together with the uniform integrability of the martingale $S$, proved in Lemma \ref{lem:moments}, this enables us to apply the criterion of Kazamaki (see Revuz, Yor \cite{revuzyor99}, p. 332). Hence we have proved the existence of a measure solution to our BSDE (\ref{BSDE}).
\hfill \epf

As a by-product of our main result, we obtain the exponential integrability of the quadratic variation of $S$.
\begin{co}\label{exponential_integrability}
Under the conditions of Theorem \ref{existence_bounded_below} we have
$$\E \exp \left( \eh\alpha^2 \int_0^T Z_s^2\, ds \right) < \infty.$$
\end{co}

\pf This follows immediately from (\ref{exponential_3}) and the lower boundedness of $\xi.$ \hfill \epf

\subsection{A quadratic BSDE with two solutions}

Let us now come back to the question of uniqueness of solutions,
and their measure solution property. Briand, Hu \cite{briandhu05}
prove the existence of solutions $(Y,Z)$ in the usual sense, given
that (\ref{eq:exp_boundedness}) is satisfied. In a setting with
more general generators the nonlinear $z$-part being bounded by
$\alpha z^2$, they provide pathwise upper and lower bounds for
$Y$, given by the known explicit solution for this generator
$(\frac{1}{2\alpha} \log \E(\exp(2\alpha \xi)|\f_t)_{t\in[0,T]}$
used above, and its negative counterpart $(-\frac{1}{2\alpha} \log
\E(\exp(-2\alpha \xi)|\f_t)_{t\in[0,T]}.$ In a more recent paper,
Briand, Hu \cite{briandhu07} also provide a uniqueness result for
the same setting, which is satisfied under the stronger
integrability hypothesis \be\label{eq:allmoments} \E(\exp(\gamma
|\xi|))< \infty \ee for all $\gamma>0$ and a convexity assumption
concerning the generator. Let us start our discussion of
uniqueness and the measure solution property by giving some
examples.

For $b>0$, let $\tau_b = \inf\{t \ge 0: W_t \le b t - 1 \}$. We
first consider a BSDE with random time horizon $\tau_b$. Let the
generator be further specified by $\alpha = \eh$. Let $\xi =
2a(b-a)\tau_b - 2a$, where $a > 0$. It will become clear along the
way why this choice of terminal variable is made. In the first
place, it is motivated by the striking simplicity of the solutions
we shall construct. We shall in fact give two explicit solutions
of the BSDE
\begin{equation} \label{counterexample}
 Y_{t \wedge \tau_b} = \xi -\int_t^{\tau_b} Z_s dW_s + \int_t^{\tau_b} \frac12 Z^2_s ds.
\end{equation}
Appropriate choices of $a$ and $b$ allow for terminal variables
that are bounded below as well as bounded above. The fact that the
time horizon is random is not crucial. Indeed, by using a time
change, any solution of Equation (\ref{counterexample}) can be
transformed into a solution of a BSDE with the same generator and
with time horizon $1$. To this end consider the time change
$\rho(t) = \frac{t}{1+t}$, $t \in [0, \infty]$, and observe that
the inverse of $\rho$ is given by $\rho^{-1}(t) = \frac{t}{1-t}$,
$t \in [0,1]$. Let $h(t) = \frac{1}{1-t}$ for all $t \in [0, 1]$.
Then the process defined by \be \label{inter} \tilde W_t =
\int_0^t h^{-1}(s) d(W_{\rho^{-1}(s)}), \qquad t \in [0,1], \ee is
a Brownian motion on $[0, 1]$. Note that $W_t = \int_0^{\rho(t)}
h(s) d \tilde W_{s}$ (and this is how we have to define $W$, if
$\tilde W$ is given). Moreover, the stopping time
\[
\hat \tau_b = \inf \left\{ t \ge 0: \int_0^{t} h(s) d \tilde W_{s}
\le \frac{t}{1-t} -1\right\}
\]
is equal to $\rho(\tau_b)$. We can now define a time changed
analogue of Equation (\ref{counterexample}) with time horizon $1$.
\begin{lem} \label{timelem}
Let $(Y_t, Z_t)$ be a solution of the BSDE (\ref{counterexample}),
and let $\hat \xi = 2a(b-a) \frac{\hat \tau_b}{1 - \hat \tau_b} -
2a$. Then $(y_t, z_t) = (Y_{\rho^{-1}(t)}, h(t) Z_{\rho^{-1}(t)})$
is a solution of the BSDE
\begin{equation} \label{counterex2}
 y_t = \hat \xi - \int_t^1 z_s d \tilde W_s + \int_t^1 \frac12 z^2_s ds.
\end{equation}
\end{lem}
\pf Since stochastic integration and continuous time changes can
be interchanged (see Proposition 1.5, Chapter V in
\cite{revuzyor99}), we have
\begin{eqnarray*}
y_t &=& Y_{\rho^{-1}(t)} = \int_0^{\rho^{-1}(t)} Z_s dW_s - \frac12 \int_0^{\rho^{-1}(t)} Z^2_s ds \\
&=& \int_0^{t} Z_{\rho^{-1}(s)} dW_{\rho^{-1}(s)} - \frac12 \int_0^{t} Z^2_{\rho^{-1}(s)} d\rho^{-1}(s) \\
&=& \int_0^{t} Z_{\rho^{-1}(s)} h(s) d \tilde W_{s} - \frac12
\int_0^{t} Z^2_{\rho^{-1}(s)} h^2(s) ds,
\end{eqnarray*}
and hence the result. \hfill \epf

Let us first assess exponential integrability properties of $\xi$.
For this, let $\gamma > 0$ be arbitrary. Then we have
\begin{eqnarray*}
\E e^{\gamma |\xi|} &=& \E e^{\gamma |2a(b-a)\tau_b - 2a|} \le
e^{2a\gamma} \E e^{\gamma 2a|b-a|\tau_b}.
\end{eqnarray*}
Define the auxiliary stopping time $$\sigma_b =  \inf\{t \ge 0:
W_t \le  t - b \}.$$ It is well known and proved by the scaling
properties of Brownian motion that the laws of $\tau_b$ and
$\frac{\sigma_b}{b^2}$ are identical (see Revuz, Yor
\cite{revuzyor99}). Moreover, the Laplace transform of $\sigma_b$
is equally well known. According to Revuz, Yor \cite{revuzyor99}
we therefore have for $\lambda > 0$ \be\label{Laplace}
E(\exp(-\lambda \tau_b)) = E(\exp(-\frac{\lambda}{b^2} \sigma_b) =
\exp(-b [\sqrt{1 + \frac{2 \lambda}{b^2}}-1]). \ee Moreover, it is
seen by analytic continuation arguments that this formula is even
valid for $\lambda \ge -\frac{b^2}{2}.$ Now choose $\lambda = -2a
|b-a| \gamma.$ Then the inequality
$$- 2a |b-a| \gamma \ge - \eh b^2$$
amounts to \be\label{example:conditions_integrability} \gamma \le
\frac{b^2}{4a |b-a|}.\ee This in turn means that we have exponential
integrability of orders bounded by $\frac{b^2}{4a |b-a| }$, in
particular we may reach arbitrarily high orders by choosing $a$ and
$b$ sufficiently close. But no combination of $a$ and $b$ allows
exponential integrability of all orders. In the light of Briand, Hu
\cite{briandhu07} this means that the entire field of pairs of
positive $a$ and $b$ promises multiple solutions, and this is
precisely what we will exhibit.\par\medskip

{\flushleft {\bf The first solution}}
\medskip

It is clear from the definition that the pair $(Y_t, Z_t)$,
defined by $Y_t = 2a W_{t \wedge \tau_b} - 2 a^2 (\tau_b \wedge t)$
and $Z = 2a 1_{[0,\tau_b]}$, is a solution of
(\ref{counterexample}). To answer the question whether this
defines a measure solution, we have to investigate

$$\E\exp \left[
\int_0^{\tau_b} \frac12 Z_s dW_s - \frac18 \int_0^{\tau_b} Z_s^2 ds
\right] = \E \exp \left[ a W_{\tau_b} - \frac{a^2}{2} \tau_b \right]
= \E(\exp(a(b-\frac{a}{2})\tau_b -a)).$$ Due to (\ref{Laplace}) we
have $$\E(\exp(a(b-\frac{a}{2})\tau_b -a)) =
\exp(-b[\sqrt{1-\frac{2}{b^2} a(b-\frac{a}{2})}-1] -a) =
\exp(-b[|1-\frac{a}{b}| - 1] -a),$$ and the latter equals 1 in case
$b\ge a$ and $\exp(2(b-a)) < 1$ in case $a>b$. This simply means
that our first solution is a measure solution of (\ref{counterex2})
provided $b\ge a$, and it fails to be one in case $a>b.$ We will
show that this first solution does not necessarily correspond to the
particular solution discussed in the beginning of the
section.\par\medskip

{\flushleft \bf The second solution}
\medskip

We show now that the BSDE (\ref{counterexample}) with the same
terminal variable as above possesses a second solution. By Lemma
\ref{timelem} there exists a second solution of (\ref{counterex2})
as well. Once this is shown, for any possible degree $\gamma$ of
exponential integrability we will have exhibited a negative random
variable satisfying $\E(\exp(\gamma |\xi|))<\infty$ for which
(\ref{counterexample}) possesses at least two solutions. This in
turn will underline that Briand, Hu's \cite{briandhu07} uniqueness
result, valid under (\ref{eq:allmoments}) cannot be improved by
much. Note that the solution we will exhibit is again of the
explicit form (\ref{solution_explicit}) encountered earlier. Let
$M_t = \E[e^\xi |\f_t]$ for all $t \ge 0$. Due to the martingale
representation property there exists a process $H$ such that $M_t =
M_0 + \int_0^t H_s dW_s$. We know that $(\ln M_{\tau_b \wedge t},
\frac{H_{\tau_b \wedge t}}{M_{\tau_b \wedge t}})$ is a solution of
(\ref{counterexample}). We will show below that
\bean
\ln M_{\tau_b \wedge
t} &=& 2b-4a + 2(b-a) W_{\tau_b\wedge t} - 2(b-a)^2 (\tau_b \wedge t),\quad \mbox{if}\quad 2a> b,\label{particular:1}\\
\ln M_{\tau_b} &=& 2a W_{\tau_b\wedge t} - 2a^2 \tau_b\wedge t,\quad \mbox{if}\quad 2a\le b.\label{particular:2}
\eean
This implies that the solution $(\ln M_{\tau_b \wedge t},
\frac{H_{\tau_b \wedge t}}{M_{\tau_b \wedge t}})$ is different from the
solution $(2a W_{\tau_b \wedge t}- 2a^2 (\tau_b \wedge t), 2a)$ obtained
above in case $2a > b$. Hence by Lemma \ref{timelem} we obtain a second solution of
(\ref{counterex2}) in this case.

First note that
\bean
M_t & =& e^{-2a} \E[ e^{2a(b-a) \tau_b}|\f_t] \\ \nonumber
&= &e^{-2a} \left( e^{2a(b-a)\tau_b} 1_{\{\tau_b \le t\}} + e^{2a(b-a)t} \E[e^{2a(b-a)[\tau_b - t]} | \f_t] 1_{\{\tau_b > t\}} \right)
\eean
Let $\sigma_b(x,t) = \inf\{s \ge 0: W_{s+t} - W_t \le b(s+t)-1-x\}$ and observe that on the set $\{\tau_b > t\}$
we have $\tau_b - t = \sigma_b(W_t, t)$.
Therefore, by using again our knowledge on the Laplace transforms of $\sigma(x,t)$ (see \cite{revuzyor99}), we get
\begin{eqnarray*}
\E[e^{2a(b-a)[\tau_b - t]} | \f_t] 1_{\{\tau_b > t\}} &=& \E[e^{2a(b-a)\sigma_b(x,t)}]\bigg|_{x = W_t} 1_{\{\tau_b > t\}} \\
&=& e^{-b(1+W_t-bt)[\sqrt{1 - \frac{4a(1-a)}{b^2}}-1]} 1_{\{\tau_b > t\}} \\
&=& e^{-b(1+W_t-bt)[|1-\frac{2a}{b}| - 1]} 1_{\{\tau_b > t\}}.
\end{eqnarray*}
Consequently,
\begin{eqnarray*}
M_t &=& e^{-2a} \left( e^{2a(b-a)\tau_b} 1_{\{\tau_b \le t\}} + e^{2a(b-a)t} e^{-b(W_t+1-bt) [|1-\frac{2a}{b}|-1]}
1_{\{\tau_b > t\}} \right)\\
&=& e^{2a((1-a)(\tau_b \wedge t) -1)} 1_{\{\tau_b \le t\}} + e^{-2(b -a)(W_t+1-bt)} 1_{\{\tau_b > t\}} .
\end{eqnarray*}
Hence in case $2a>b$
\begin{eqnarray*}
\ln M_{\tau_b \wedge t} &=& 2a((b-a)(\tau_b \wedge t) -1) - 2(a -b)(W_{\tau_b \wedge t}+1-(\tau_b \wedge t)) \\
&=& -4a + 2b + [-2b +4a - 2a^2] (\tau_b \wedge t) - 2(a-b) W_{\tau_b \wedge t} \\
 &=& 2b - 4a + 2(b-a) W_{\tau_b \wedge t}- 2[b-a]^2 (\tau_b \wedge t) .
\end{eqnarray*}

This confirms the first equation (\ref{particular:1}). Let finally $2a\le b.$
Then we have

\begin{eqnarray*}
M_t &=& e^{-2a} \left( e^{2a(b-a)\tau_b} 1_{\{\tau_b \le t\}} + e^{2a(b-a)t} e^{2a (W_t+1-bt)}
1_{\{\tau_b > t\}} \right)\\
&=& e^{2a((b-a)(\tau_b \wedge t) +2a(W_{\tau_b\wedge t}+1-b\tau_b\wedge t)}\\
&=& e^{2a W_{\tau_b\wedge t} - 2 a^2 \tau_b\wedge t}.
\end{eqnarray*}
Hence in this case
\begin{eqnarray*}
\ln M_{\tau_b \wedge t} &=& 2a W_{\tau_b \wedge t} - 2a^2 \tau_b \wedge t.
\end{eqnarray*}
Note that in case $2a\le b$ we recover the solution already obtained as the first solution.

Let us finally show that this second solution is in fact a measure solution for any possible combination of parameters.

\begin{lem}
$(\ln M_{\tau_b \wedge t}, \frac{H_{\tau_b \wedge t}}{M_{\tau_b \wedge t}})$ can be extended to a measure solution of (\ref{counterexample}),
hence provides a measure solution of (\ref{counterex2}).
\end{lem}
\pf
For the first solution in case $a\le b$, which is identical to the one considered in case $2a\le b$, we have already established
the measure solution property. Let us therefore consider the case $2a>b$.
Note that for all $t$, $M_{t \wedge \tau_b} = e^{2b -4a} + \int_0^{t \wedge \tau_b} H_s dW_s$. It\^o's formula
applied to $e^{2(b-a) W_{\tau_b \wedge t}- 2[b-a]^2 (\tau_b \wedge t) }$ yields
\[ H_{s \wedge \tau_b} = 2(b-a) e^{2(b-a) W_{\tau_b \wedge t}- 2[b-a]^2 (\tau_b \wedge t) }. \]
As a consequence, we have
\[Z_{s\wedge \tau_b} = \frac{H_{s\wedge \tau_b}}{M_{s \wedge \tau_b}} = 2(b-a) 1_{[0, \tau_b]}(s), \]
and therefore
\begin{eqnarray*}
\mathcal E \Big(\frac12 \int Z dW\Big)_{\tau_b} &=& e^{(b-a) W_{\tau_b} -\frac12 (b-a)^2 \tau_b} \\
&=&  e^{(b-a) (b \tau_b -1) -\frac12 (b-a)^2 \tau_b} \\
&=& e^{(a-b)} e^{\frac12 (b-a)(b+a) \tau_b}.
\end{eqnarray*}
Again the explicit representation of the Laplace transform in (\ref{Laplace}) gives
\[ \E \mathcal E \Big(\frac12 \int Z dW\Big)_{\tau_b} = e^{(a-b)} \E e^{- \frac12 (b-a)(b+a) \tau_b}
= e^{(a-b)} e^{- b(\sqrt{1 - (1-\frac{a^2}{b^2})} -1)} = 1. \]
This implies the claimed result that our second solution $(\ln M_{\tau_b \wedge t}, \frac{H_{\tau_b \wedge t}}{M_{_b \wedge t}})$
is a measure solution of (\ref{counterexample}).
\hfill \epf

{\flushleft \bf Remarks:}

1. We can summarize the findings of our investigations of the examples by stating that there are
three basic scenarios: a) for $b\ge 2a$ we obtained one solution which is a measure solution at the same time;
b) in the range $2a>b\ge a$ we found two different solutions both of which are measure solutions;
c) if $a>b$ we finally encountered two solutions one of which is a measure solution, while the other one
is not.\par\smallskip

2. Note that our examples exhibiting solutions of
(\ref{counterexample}) that are not measure solutions are all for
negative terminal variables $\xi.$ Positive terminal variables arise
in scenarios a) or b), and therefore only produce multiple measure
solutions.\par\medskip

{\flushleft \bf A continuum of solutions}
\medskip

Let us now combine the first and second solutions to obtain a
continuum of solutions of our BSDE (\ref{counterexample}). To do
this, we have to consider a still somewhat more general class of
stopping times. For $c\in\r$, let
$$\rho_c = \inf \{t\ge 0: W_t \le t - c\}.$$
We investigate the terminal variables
$$\xi = 2a (a-1) \rho_c + d$$
with further constants $a\not=0, d\in\r.$ Note first that the
integrability properties of $\xi$ are the same as those obtained
before for $b=1.$ According to the preceding paragraphs, our BSDE
(\ref{counterexample}) possesses the following two solutions
\be\label{solution:1}
Z^1 = 2a 1_{[0,\rho_c]},\quad Y^1 = d_1 + 2a W_{\rho_c \wedge \cdot}- 2 a^2 \rho_c\wedge \cdot,\ee
\be\label{solution:2}
Z^2 = 2(1-a) 1_{[0,\rho_c]},\quad Y^2 = d_2 + 2(1-a) W_{\rho_c \wedge \cdot}- 2 (1-a)^2 \rho_c \wedge \cdot,\ee
with $d_1 = -2ac$ resp. $d_2 = -2(a-1)c$.
Let us now take $c=1$ and combine the two solutions to obtain a continuum of new ones. To do this, we start with the
equation
$$\rho_1 = \rho_c + \rho_{1-c}\circ \theta_{\rho_c},$$
where $theta_t$ is the shift on Wiener space defined by
$$\theta_t(\omega) = W_{t+\cdot}(\omega) - W_t(\omega),$$
and $c\in ]0,1[.$ It describes the first time to reach the line with
slope 1 that cuts the vertical at level $-1$, by decomposition with the intermediate time to reach the line with slope 1
cutting the vertical at $-c$. We mix the two solutions on the two resulting
stochastic intervals, more precisely we put for $c\in]0,1[, l\in\r$
\bean\label{solution:mixed}
Z^c &=& 2a 1_{[0,\rho_c]} + 2(1-a)1_{[\rho_c, \rho_1]},\\\nonumber
Y^c &=& l+ 2a W_{\rho_c\wedge\cdot} - 2a^2 \rho_c\wedge\cdot + 2(1-a)[W_{\rho_1\wedge\cdot}-W_{\rho_c\wedge\cdot}] - 2(1-a)^2
[\rho_1\wedge\cdot  - \rho_c\wedge\cdot].
\eean
Since we have
\bea Y^c_{\rho_1} &=& l+ 2a W_{\rho_c} - 2a^2 \rho_c + 2(1-a)[W_{\rho_1}-W_{\rho_c}] - 2(1-a)^2
[\rho_1  - \rho_c]\\ &=& l + 2a(1-a) \rho_1 - 2ac - 2(1-a)(1-c),\eea
we have to set
$$l-2ac -2(1-a)(1-c) = d$$
in order to obtain a solution of (\ref{counterexample}) with $c=1.$
According to the treatment of the first and second solution, the constructed mixture is
a measure solution if and only if both components of the mixture are. This is the case for
$2a(1-a) > 0$, whereas for $2a(1-a)<0$ we obtain a continuum of solutions that are no measure solutions.\par\bigskip

{\flushleft \bf Remarks:}

1. This time, we may summarize our results by saying that there are two scenarios:
a) for $2a(1-a)>0$ there is a continuum of measure solutions of (\ref{counterexample}), while for
$2a(1-a)<0$ a continuum of non measure solutions is obtained.

2. Note that the initial conditions of our solutions continuum vary in a convex way
between $-2a$ and $-2(1-a)$ as $c$ varies in $]0,1[,$ spanning the whole interval.\par\bigskip

We shall now point out that the measure solution property of the second solution in case $a>b$ exhibited
in the example above is not a coincidence. In fact,
it will turn out that also for negative exponentially integrable $\xi$, solutions given by
(\ref{solution_explicit}) provide measure solutions. To prove this, we will reverse the sign of $\xi$ by
looking at our BSDE from the perspective of an equivalent measure.

\subsection{Exponentially integrable upper bounded terminal variable} \label{subsection:upperexponential}

Sticking with the positivity of $\alpha$ in the generator
$$f(s,z) = \alpha z^2,\quad s\in[0,T], z\in\mathbb{R}$$
we shall now consider terminal variables $\xi$ that fulfill the
exponential integrability condition (\ref{eq:exp_boundedness}),
but are bounded above by a constant. Again, by a constant shift of
the solution component $Y$, we can assume that the upper bound is
0, i.e. $\xi \le 0.$
So fix a non-positive terminal variable $\xi$ satisfying (\ref{eq:exp_boundedness}) for some $\gamma > 2 \alpha,$ and denote by $(Y,Z)$ the pair of processes given by the explicit representation of (\ref{solution_explicit}) solving our BSDE according to Briand, Hu \cite{briandhu05}.
With respect to the following probability measure, $\xi$ will effectively change its sign, so that we can hook up to the previous discussion. Recall $S = \int_0^\cdot Z_s d W_s.$

\begin{lem}\label{lem:signchange}
Let $V = \exp(2\alpha S - 2 \alpha^2 \langle S\rangle).$ Then $V$ is a martingale of class (D), and consequently
$$R = V_T \cdot \prb$$
is a probability measure equivalent to $\prb.$ Moreover,
$$W^R = W - 2 \alpha \int_0^\cdot Z_s ds$$
is a Brownian motion under $R$.
\end{lem}
\pf By (\ref{BSDE}), we may write
$$2 \alpha [Y-Y_0] = 2 \alpha S - 2 \alpha^2 \langle S\rangle,$$
hence also
$$2 \alpha [\xi-Y_0] = 2 \alpha S_T - 2 \alpha^2 \langle S\rangle_T.$$
According to Briand, Hu \cite{briandhu05}, Theorem 2, there exists $\delta > 2\alpha$ such that
\be\label{class_D}
\E(\sup_{t\in[0,T]} \exp(\delta |Y_t|)) < \infty,
\ee
and therefore $\beta > 1$ with the property
\be\label{class_D2}
\E(\sup_{t\in[0,T]} V_t^{\beta}) < \infty.
\ee
This clearly implies that $V$ is a martingale of class (D), and consequently $R$ is a probability measure. Finally, Girsanov's theorem implies that $W^R$ is a Brownian motion under $R.$\hfill \epf

Now consider our BSDE under the perspective of the measure $R$ with respect to the Brownian motion $W^R.$ We may write
\be\label{switch}
Y = \xi - \int_\cdot^T Z_s d W_s + \alpha \int_\cdot^T Z_s^2 ds = \xi - \int_\cdot^T Z_s dW_s^R - \alpha \int_\cdot^T Z_s^2 ds.
\ee
But this just means that by switching signs in $(Y,Z)$, we may return, under the new measure $R$, to our old BSDE with $\xi$ replaced with $-\xi.$ So our measure change puts us back into the framework of the previous subsection, and we may resume our arguments there by setting
$$S^R = -\int_0^\cdot Z_s d W^R_s.$$
We need an analogue of Lemma \ref{lem:moments} to guarantee that $R$ is a uniformly integrable martingale.

\begin{lem}\label{lem:moments_switch}
For any $p\ge 1$ we have
$$\E^R \left( \left[ \int_0^T Z_s^2 ds\right]^p \right) < \infty.$$
In particular, $S^R$ is a uniformly integrable martingale under $R.$
\end{lem}
\pf By definition of $R$, we have for any $p>1$
$$\E^R \left( \left[ \int_0^T Z_s^2 ds\right]^p\right) = \E \left( \exp(2\alpha [\xi - Y_0]) \left[ \int_0^T Z_s^2 ds \right]^p \right).$$
Now since $\xi\le 0,$ the density $\exp(2\alpha [\xi - Y_0])$ is bounded above. Therefore the asserted moment finiteness follows from
Lemma \ref{lem:moments}. \hfill \epf

We are in a position to prove the main result of this subsection.

\begin{thm}\label{existence_bounded_above}
Assume that that $f$ satisfies $f(s,z) = \alpha z^2,$ $ z\in\r, s\in[0,T]$, and that $\xi$ is bounded above and satisfies (\ref{eq:exp_boundedness}). Then there is a measure solution of (\ref{BSDE}) with a measure $\qprb$ that is equivalent to $\prb$.
\end{thm}
\pf We may assume $\xi\le 0.$ Let us first show, in analogy to the proof of Theorem \ref{existence_bounded_below}, that
$$V^R = \exp(\alpha S^R - \eh \alpha^2 \langle S^R\rangle)$$
is a uniformly integrable martingale under $R$, using Kazamaki's criterion. For this purpose, let
$$\tau_n^R = \inf\{ t\ge 0: \langle S^R\rangle_t \ge n\}\wedge T,\quad n\in\mathbb{N}.$$
Then, due to $\langle S\rangle = \langle S^R\rangle,$ we deduce for all $n\in\mathbb{N}$
that $\tau_n = \tau_n^R.$ Since $\tau_n^R \to T$ as $n\to\infty$, even with $\tau_n^R = T$ for all but finitely many $n$, Fatou's lemma allows to deduce
\be\label{fatou_switch}
\E^R(V_T) \le \liminf_{n\to\infty} \E^R(V^R_{\tau_n^R}) \le 1.
\ee
Moreover, by the form of the BSDE translated to $W^R$ under $R$,
\bea
\alpha S^R - \eh \alpha^2 \langle S^R\rangle &=& \alpha [-\int_0^\cdot Z_s d W^R_s - \eh \alpha \int_0^\cdot Z_s^2 ds]\\
&=& \alpha[-\int_0^\cdot Z_s d W^R_s - \alpha \int_0^\cdot Z_s^2 ds] + \eh \alpha^2 \int_0^\cdot Z_s^2 ds\\
&=& \alpha [-Y+Y_0] + \eh \alpha^2 \int_0^\cdot Z_s^2 ds.
\eea
Using the negativity of $\xi$ and the identity just derived, we get the integrability property
\be\label{exponential_3_switch}
\E^R\exp \left[ \eh \alpha (-\xi + Y_0) + \eh \alpha^2 \int_0^T Z_s^2\, ds \right] < \infty.
\ee
Using this and the positivity of the terminal variable $\xi$, we can now obtain the exponential integrability property
\be\label{exponential_4} \E\exp \left[ \eh \alpha (\xi - Y_0)  + \eh \alpha^2 \int_0^T Z_s^2\, ds \right] < \infty. \ee
Again, we may now use (\ref{exponential_3_switch}) together with (\ref{fatou_switch}) to prove the exponential integrability of $\eh \alpha S^R_T$. In fact, from the BSDE viewed with $W^R$ under $R$ we have
$$\eh \alpha S_T^R =  \eh\alpha (-\xi + Y_0) + \eh \alpha^2 \int_0^T Z_s^2\, ds.$$
Hence we obtain
\be\label{kazamaki2} \E^R \exp \left( \eh \alpha S^R_T \right) < \infty. \ee
Now appeal to the uniform integrability of the martingale $S^R$ under $R$, proved in Lemma \ref{lem:moments_switch}, to see that the criterion of Kazamaki (see Revuz, Yor \cite{revuzyor99}, p. 332) may be applied. Hence $V_R$ is a uniformly integrable martingale under $R$.

We have to show that this implies uniform integrability of
$$V = \exp(\alpha S - \eh \alpha^2 \langle S\rangle)$$
under $\prb.$ To see this, note that
\bea
\exp(\alpha S - \eh\alpha^2 \langle S\rangle) &=& \exp(2\alpha S - 2\alpha^2 \langle S\rangle)\cdot \exp(-\alpha S +\frac{3}{2}\alpha^2 \langle S\rangle)\\
&=& \exp(2\alpha S - 2\alpha^2 \langle S\rangle)\cdot \exp(\alpha S^R - \eh\alpha^2 \langle S^R\rangle).
\eea
Hence for $n\in\mathbb{N}$
\be\label{shiryaev_swith}
\E(V_{\tau_n} 1_{\{\tau_n<T\}}) = \E^R(V_{\tau_n^R}^R 1_{\{\tau_n^R < T\}}),
\ee
and the latter expression tends to 0 as $n\to \infty$ by the first part of the proof. Hence the uniform integrability of $V$ under $\prb$ follows from the explosion criterion (\ref{explosion}) already used earlier. This completes the proof. \hfill \epf

{\flushleft \bf Remark:} The results of the preceding two subsections clearly call for
similar ones for our BSDE with exponentially integrable terminal
variable that are not bounded. Due to the nonlinearity of the
generator of the BSDE, it seems impossible to derive such properties
by combining the results of Theorems \ref{existence_bounded_below}
and \ref{existence_bounded_above}.

\section{The existence of measure solutions in the Lip\-schitz case}

We shall now construct measure solutions from first principles. In
particular, we shall not assume any knowledge about strong
solutions. We shall give a complete self-contained construction for
measure solutions with Lipschitz continuous generator for which the
Lipschitz constant may be time dependent. Our construction provides
the measure solution along an algorithm which just iterates the
procedures of projecting the terminal variable by a given measure.
The martingale representation theorem with respect to the measure
$\mathbb{Q}^n$ in step $n$ will produce a control process $Z^n$
which is then plugged into the generator of the BSDE. The resulting
drift is taken off by applying Girsanov's theorem which produces a
new measure $\mathbb{Q}^{n+1}$ with which we continue along the
lines just sketched in step $n+1.$ The sequence
$(\mathbb{Q}^n)_{n\in{\bf N}}$ thus produced has to be shown to
possess at least an accumulation point in the weak topology. This is
seen by a simple argument using the Lipschitz and boundedness
properties. The extension to a continuous or quadratic generator and
bounded terminal condition is straightforward from this perspective,
and uses monotone approximations following the scheme in
\cite{lepeltiersanmartin98}. But this result is already contained in
the results of \cite{kobylanski00} and Theorem
\ref{existence_bounded_liptser}. Hence we do not write the details
here. The extension of our intrinsic construction of measure
solutions to unbounded terminal conditions is left for future
research.

In order to obtain a self-contained theory that is not using any
knowledge on classical solutions, we first construct measure
solutions in a setting for which they have been studied mostly: for
generators that increase at most linearly and possess Lipschitz
properties with time dependent and random Lipschitz constants. More
formally, in this section we consider the following class of
generators. Let
$$f:\Omega \times [0,T] \times \R \to \R$$ satisfy the
{\bf Assumption (H2)}: for some $\gamma \ge 1$ and some non negative process $\phi$
\begin{enumerate}
\item $\xi \in L^{\gamma}(\Omega)$;
\item $f(s,z) = f(\cdot,s,z)$ is adapted for any $z \in \R$;
\item $\displaystyle \E \left( \int_0^T |f(s,0)|^{\gamma} ds \right) < \infty$;
\item the set $\left\{ s \in [0,T], \ f(s,.) \ \mbox{is not continuous} \right\}$ is of Lebesgue measure zero;
\item $|f(s,z) - f(s,z')| \le \phi_s |z-z'|$ for all $s \in [0,T], \  (z,z')\in \R^2$.
\end{enumerate}
We shall assume in the following that $f(s,0) = 0$ for all $s \in [0,T]$. This can be done without loss of generality, since we may replace $\xi$ with the $\gamma$-integrable random variable
$$\tilde{\xi} = \xi + \int_0^T f(s,0) ds.$$

Now we define the function $g : \Omega \times [0,T] \times \R \to \R$ by the requirement that for all $s \in [0,T], \ z \in \R$:
\begin{eqnarray*}
g(s,z) & = & \frac{f(s,z)}{z},\ \mbox{if} \ z \neq 0, \\
& = & 0,\ \mbox{if} \ z = 0.
\end{eqnarray*}
Therefore we have defined the function $g$ with values in $\R$ and $g$ is bounded by the process $\phi$.

The process $\phi$ verifies either
\begin{equation} \label{hyp1}
\exists \kappa > 1, \ \E \left[ \exp \left( \frac{\kappa}{2}
\int_{0}^{T} \phi_{r}^{2} dr \right) \right] < + \infty
\end{equation}
or
\begin{equation} \label{hyp2}
\mbox{the martingale} \ \left( L_{t} = \int_{0}^{t} \phi_{r} dW_{r}\right)_{t \in [0,T]} \ \mbox{is BMO}.
\end{equation}
 We denote by $\|L\|$ the $BMO_{2}$-norm of $L$. From Theorem 2.2 in \cite{Kaz}, (\ref{hyp2}) implies (\ref{hyp1}), with $1/\kappa = 2 \|L\|^2$. Remark that (\ref{hyp1}) is a stronger Novikov condition. From these assumptions (see \cite{Kaz}, theorem 2.3), we know that for $0 \leq t \leq T$,
$$\expo( \phi W)_{t} = \exp \left( \int_{0}^{t} \phi_{r} dW_{r} - \frac{1}{2} \int_{0}^{t} \phi_{r}^{2} dr \right)$$
is a uniformly integrable martingale.

We define the process $\Phi$ by
$$\forall t \in [0,T], \quad \Phi_t = \int_0^t \phi_s^2 ds,$$
and {\bf Assumption (H3)} holds: there exists two constants $\alpha>
\Psi$ and $\delta> \Psi$ such that
\begin{equation} \label{hyp3}
\E(e^{\alpha \Phi_T} |\xi|^{\delta}) < + \infty.
\end{equation}
The constant $\Psi>1$ is given for (\ref{hyp1}) by:
$$\Psi(\kappa) = \Psi_{\ref{hyp1}}(\kappa) = 1+ 4 \frac{\sqrt{\kappa} }{(\sqrt{\kappa}-1)^2}= \left( 1 + \frac{2\sqrt{\kappa}+1}{\kappa} \right)  \frac{\kappa }{(\sqrt{\kappa}-1)^2},
$$
and for (\ref{hyp2}) by:
$$\Psi(\|L\|) = \Psi_{\ref{hyp2}}(\| L \|) = \left( 1 + \frac{\|L\|}{2} \right)  \frac{\theta^{-1} \left( \| L \| \right)}{\theta^{-1} \left( \| L \| \right) -1}.$$
The function $\theta : ]1,+\infty[ \to \R_{+}^{*}$ is the continuous decreasing function given by
$$\forall q \in ]1,+\infty[, \ \theta(q) = \left\{ 1 + \frac{1}{q^{2}} \ln \frac{2q-1}{2(q-1)} \right\}^{\frac{1}{2}} -1.$$
We can check that $\Psi_{\ref{hyp2}} : ]0,+\infty[ \to ]1,+\infty[$ is an increasing function such that $\Psi(0)=1$ and $\Psi(\infty) = \infty$.

\begin{rem}
If $f$ is a Lipschitz function:
$$|f(t,z)-f(t,z')| \leq K |z-z'|,$$
then $\phi$ is the constant $K$. Then (\ref{hyp1}) is satisfied for all $\kappa > 1$, and (\ref{hyp3}) holds if $\gamma > 1$.
\end{rem}

The solution algorithm for our BSDE (\ref{BSDE})
$$Y_t= \xi + \int_t^T f(s,Z_s) ds - \int_t^T Z_s dW_s$$
is based on a recursively defined change of measure. Let
$\mathbb{Q}^0 = \prb,$ and $W^0 = W$, the coordinate process which
is a Wiener process under $\mathbb{Q}^0$. Set
$$Y^1 = \E (\xi|\tri_\cdot) = \E (\xi) + \int_0^\cdot Z^1_s d W_s^0,$$
and
$$\mathbb{Q}^1 = \exp \left( \int_0^T g(s, Z^1_s) d W_s - \frac{1}{2} \int_0^T g(s, Z_s^1)^2 ds \right) \cdot \prb = R^1_T \cdot \prb.$$
Then
$$W^1 = W - \int_0^\cdot g(s, Z_s^1) ds$$
is a Wiener process under $\mathbb{Q}^1$. Indeed under (\ref{hyp1}),
the Novikov condition is satisfied, and under (\ref{hyp2}), the
martingale
$$M^1_t = \int_0^t g(s,Z^1_s) dW_s$$
is BMO. Now since $(\mathbb{Q}^1, \mathbb{Q}^0)$ is a Girsanov pair,
it is well known that the predictable representation property is
inherited from the Brownian motion $W^0$ to the Brownian motion
$W^1$. See for example Revuz, Yor \cite{revuzyor99}, p. 335. Hence
there exists a pair $(Y^2, Z^2)$ of processes such that for all
$t\in[0,T]$
$$Y^2_t = \E^{\mathbb{Q}^1} (\xi | \tri_t) = \E^{\mathbb{Q}^1}(\xi) + \int_0^t Z^2_s d W_s^1.$$

Assume that $\mathbb{Q}^n$ is recursively defined, along with the
Brownian motion
$$W^n = W - \int_0^\cdot g(s, Z^n_s) ds$$
under $\mathbb{Q}^n$. Then Revuz, Yor \cite{revuzyor99} may be
applied to obtain two processes $(Y^{n+1}, Z^{n+1})$ such that
$$Y^{n+1} = \E^{\mathbb{Q}^n}(\xi|\tri_\cdot) = \E^n(\xi|\tri_\cdot) = \E^n (\xi) + \int_0^\cdot Z_s^{n+1} d W_s^n.$$
Now set
$$\mathbb{Q}^{n+1} = \exp \left[ \int_0^T g(s, Z^{n+1}_s) d W_s -
\int_0^T g(s, Z_s^{n+1})^2 ds \right] \cdot \prb = R^{n+1}_T \cdot
\prb$$ to complete the recursion step. Then from our assumptions on
$\phi$, and from the boundedness of $g$ by $\phi$, the sequence of
probability measures $(\mathbb{Q}^n)_{n\in\N}$ is well defined and
consists of measures equivalent with $P$. It is not hard to show
tightness for this sequence.

\begin{pr} \label{tightness_lipschitz}
Under (\ref{hyp1}) or (\ref{hyp2}), the sequence
$(\mathbb{Q}^n)_{n\in\N}$ is tight.
\end{pr}

\pf In this proof, $\E^n$ denotes the expectation under
$\mathbb{Q}^n$. For $0 \leq s \leq t \leq T$, $n \in \N$, we have,
recalling that $W$ is the coordinate process on the canonical space:
\begin{eqnarray} \nonumber
\E^{n} \left( |W_{t}-W_{s}|^{4} \right) & \leq & \E^{n} \left(|W^{n}_{t}-W^{n}_{s} + \int_{s}^{t} g(u,Z^{n}_{u}) du|^{4} \right) \\ \nonumber
& \leq & C \left[ \E^{n} \left(|W^{n}_{t}-W^{n}_{s}|^{4} \right) + \E^{n} \left( |\int_{s}^{t} g(u,Z^{n}_{u}) du|^{4} \right) \right] \\ \nonumber
& \leq & C |t-s|^{2} + C |t-s|^{2} \E^{n} \left( \int_{s}^{t} g(u,Z^{n}_{u})^{2} du \right)^{2}  \\ \nonumber
& \leq & C |t-s|^{2} + C |t-s|^{2} \E^{n} \left( \int_{s}^{t} \phi_{u}^{2} du \right)^{2} \\ \nonumber
& \leq & C |t-s|^{2} + C |t-s|^{2} \E \left[ \left( \int_{s}^{t} \phi_{u}^{2} du \right)^{2} R^{n}_T \right] \\ \label{210206-1}
& \leq & C |t-s|^{2} \left\{ 1 + \left[ \E \left( \int_{s}^{t} \phi_{u}^{2} du \right)^{2p} \right]^{1/p} \left[ \E (R^{n}_T)^{q} \right]^{1/q} \right\},
\end{eqnarray}
from the H\"{o}lder inequality with $p>1$ and $p^{-1}+q^{-1}=1$.

Suppose that $\phi$ satisfies the assumption (\ref{hyp1}). From the Novikov condition applied to the martingale
$$M^{n}_{t} = \int_{0}^{t} g(u,Z^{n}_{u}) dW_{u},$$
we know that $\expo (M^{n})$ is a uniformly integrable martingale under $\prb$. Moreover if $C \leq \kappa$
$$\E \left[ \exp \left( \frac{\sqrt{C}}{2} M^{n}_{T} \right) \right] \leq \E \left[ \exp \left( \frac{C}{2} \langle M^{n} \rangle_{T} \right) \right]^{1/2} \leq \E \left[ \exp \left( \frac{\kappa}{2} \langle L \rangle_{T} \right) \right]^{1/2} < + \infty.$$
From Theorem 1.5 in \cite{Kaz}, we deduce that if $p > p^*$ with
$$\frac{\sqrt{p^*}}{\sqrt{p^*} -1 } =  \sqrt{\kappa} \Longleftrightarrow p^*=\frac{\kappa}{(\sqrt{\kappa}-1)^2},$$
then for $q < q^*$
\begin{equation} \label{210206-2}
\E \left[ \expo \left( g(.,Z^{n}) W \right)_{T}^{q} \right] = \E (R^{n}_T)^{q} \leq C.
\end{equation}

Now if $\phi$ verifies the assumption (\ref{hyp2}), the martingale $M^{n}$ is also BMO, and the BMO-norm of $M^{n}$ is smaller than the BMO-norm of $L$. Therefore from Theorem 3.1 in \cite{Kaz} (or more precisely from the proof of this result), we deduce that there exists $q > 1$ and $C$ s.t.
\begin{equation} \label{210206-3}
\E \left[ \expo \left( g(.,Y^{n},Z^{n}) W \right)_{T}^{q} \right]= \E (R^{n}_T)^{q} \leq C.
\end{equation}
The constant $q$ must satisfy the following inequality: $q < q^*$ with
$$\|L\| = \theta(q^*) \Longleftrightarrow q^* = \theta^{-1}(\|L\|) \Longleftrightarrow p^* = \frac{\theta^{-1} \left( \| L \| \right)}{\theta^{-1} \left( \| L \| \right) -1}.$$
Moreover, from the John-Nirenberg inequality (see \cite{Kaz}, Theorem 2.2):
$$\E \left[ \exp \left( \frac{1}{4 \| L \|^{2}_{BMO_{2}}} \int_{0}^{T} \phi_{u}^{2} du \right) \right] \leq 2 \Longrightarrow \E \left( \int_{s}^{t} \phi_{u}^{2} du \right)^{2p} < + \infty.$$
Finally from (\ref{210206-1})
$$\E^{n} \left( |W_{t}-W_{s}|^{4} \right)  \leq C |t-s|^{2}.$$
Hence by a well known criterion (see for example Kallenberg \cite{kallenberg97}, p. 261), tightness follows.
\epf

In a second step, we shall now establish the boundedness in $L^2$ of the control sequence $(Z^n)_{n\in\N}$ obtained by the algorithm. Before let us give some estimates.

\begin{lem} \label{bound2_lem}
If $\delta > \Psi$ and (\ref{hyp3}) holds, there exist two constants $\beta > 0$ and $p > 1$ such that
\begin{equation} \label{hyp3bis}
\forall n \in \N, \quad \E^{n-1} \left(e^{\beta \Phi_T} |\xi|^{p} \right) < + \infty.
\end{equation}
\end{lem}
\pf In the proof of Proposition \ref{tightness_lipschitz}, we
already see that there exists $q^* > 1$ such that for every $1 < r
< q^*$, and for every $n$, $\E \left( R^{n-1}_T \right)^r \leq C_r
< + \infty$. Thus using H\"older's inequality
$$\E^{n-1} \left(e^{\beta \Phi_T} |\xi|^{p} \right) \leq \left[ \E \left(e^{s\beta \Phi_T} |\xi|^{sp} \right) \right]^{1/s} \times \left[ \E \left( R^{n-1}_T \right)^r \right]^{1/r} \leq C_r \left[ \E \left(e^{s\beta \Phi_T} |\xi|^{sp} \right) \right]^{1/s}.$$
From (\ref{hyp3}), $\delta > \Psi$ implies that $\delta > p^* = ( 1 - 1/q^*)^{-1}$. Hence for $r < q^*$, $\delta /s > 1$ and we can find $p > 1$ such that $sp < \delta$. Then choosing $\beta$ sufficiently small, $s \beta < \alpha$ and the conclusion follows.
\epf

From Lemma \ref{bound2_lem}, we deduce:
\begin{lem} \label{bound1_lem}
There exists a constant $C$ such that for every $n \in \N$,
$$\E^{n-1} \left[ \sup_{t \in [0,T]}(e^{\beta \Phi_t} |Y^n_t|^{p}) + \left( \int_0^T e^{\beta \Phi_t} |Z^n_t|^2 dt \right)^{p/2} \right] \leq C \E^{n-1} \left[ \exp \left( \beta \Phi_T  \max \left(\frac{p}{2} , 1 \right) \right) |\xi|^{p} \right].$$
\end{lem}
\pf
Recall that for every $n$, $Y^n_t = \E^{n-1} (\xi | \tri_t) = \xi - \int_t^T Z^n_s dW^{n-1}_s$. Therefore
$$e^{(\beta/p) \Phi_t} |Y^n_t| \leq  \E^{n-1} (e^{(\beta/p) \Phi_t} |\xi| | \tri_t) \leq   \E^{n-1} (e^{(\beta/p) \Phi_T} |\xi| | \tri_t) .$$
Using Doob's inequality we deduce
$$\E^{n-1} \sup_{t \in [0,T]}(e^{\beta \Phi_t} |Y^n_t|^{p}) \leq C_{p}  \E^{n-1} (e^{\beta \Phi_T} |\xi|^{p}) .$$
Now we have
$$\int_t^T e^{\beta \Phi_s/2} Z^n_s dW^{n-1}_s = e^{\beta \Phi_T/2} \xi - e^{\beta \Phi_t/2 } Y^n_t - (\beta/2) \int_t^T e^{\beta \Phi_s/2} Y^n_s \phi_s^2 ds.$$
Using Burkholder--Davis--Gundy's inequality and the previous estimate on $Y^n$
$$\E^{n-1} \left[  \left( \int_0^T e^{\beta \Phi_t} |Z^n_t|^2 dt \right)^{p/2} \right] \leq C \E^{n-1} (e^{\beta \Phi_T p/2} |\xi|^{p}).$$
\epf

\begin{pr}\label{boundedness}
Under Assumption (H2), if $\delta > \Psi$ and if (\ref{hyp3}) holds, there exists $\beta > 0$ and $p > 1$ such that
$$\E \left( \sup_{t \in [0,T]} e^{\beta \Phi_t} (Y^{n}_{t})^{p} \right) \quad \mbox{and}\quad \E \left[ \left( \int_{0}^{T} e^{\beta \Phi_s} (Z^{n}_{s})^{2} ds \right)^{\frac{p}{2}} \right]$$
are bounded sequences.
\end{pr}

\pf
We give just the proof for the sequence $\displaystyle \E \left[ \left( \int_{0}^{T} e^{\beta \Phi_s} (Z^{n}_{s})^{2} ds \right)^{\frac{p}{2}} \right]$. For the other sequence, the sketch is the same. Denote for $n\in\N$
$$R^n = R^n_T = \exp \left( \int_0^T g(s, Z^{n}_s) d W_s - \frac{1}{2} \int_0^T g(s, Z_s^{n})^2 ds \right).$$
Then for $p>1$ and $\eps > 0$
\begin{eqnarray*}
\E \left( \int_{0}^{T} e^{\beta \Phi_s} (Z^{n}_{s})^{2} ds \right)^{\frac{p}{2}} & = & \E \left[ \left( \int_{0}^{T} e^{\beta \Phi_s} (Z^{n}_{s})^{2} ds \right)^{\frac{p}{2}} (R^{n-1})^{\frac{1}{1+\eps}} (R^{n-1})^{-\frac{1}{1+\eps}} \right] \\
& \leq & \left[ \E  \left( \int_{0}^{T} e^{\beta \Phi_s} (Z^{n}_{s})^{2} ds \right)^{\frac{p(1+\eps)}{2}}  R^{n-1} \right]^{\frac{1}{1+\eps}} \left[ \E (R^{n-1})^{-\frac{1}{\eps}} \right]^{\frac{\eps}{1+\eps}} \\
& = & \left[ \E^{n-1}  \left( \int_{0}^{T} e^{\beta \Phi_s} (Z^{n}_{s})^{2} ds \right)^{\frac{p(1+\eps)}{2}} \right]^{\frac{1}{1+\eps}} \left[ \E (R^{n-1})^{-\frac{1}{\eps}} \right]^{\frac{\eps}{1+\eps}}.
\end{eqnarray*}
With Lemmas  \ref{bound1_lem} and  \ref{bound2_lem} we obtain:
$$\E^{n-1} \left( \int_{0}^{T} e^{\beta \Phi_s} (Z^{n}_{s})^{2} ds \right)^{\frac{p(1+\eps)}{2}} \leq C \E^{n-1} \left(e^{\beta \Phi_T\frac{p(1+\eps)}{2}} \left|  \xi \right|^{p(1+\eps)} \right). $$
Thus for some $\eta > 0$
\begin{eqnarray} \label{210206-4} \nonumber
&& \E \left( \int_{0}^{T} e^{\beta \Phi_s} (Z^{n}_{s})^{2} ds \right)^{\frac{p}{2}} \leq C \left[ \E \left( e^{\beta \Phi_T \frac{p(1+\eps)}{2}}\left| \xi \right|^{p(1+\eps)} R^{n-1} \right) \right]^{\frac{1}{1+\eps}} \left[ \E (R^{n-1})^{-\frac{1}{\eps}}  \right]^{\frac{\eps}{1+\eps}} \\
&& \leq C \left\{ \E e^{\beta \frac{p(1+\eps)}{2}(1+\eta) \Phi_T}\left|  \xi \right|^{p(1+\eps)(1+\eta)} \right\}^{\frac{1}{(1+\eps)(1+\eta)}} \left\{ \E (R^{n-1})^{\frac{1+\eta}{\eta}} \right\}^{\frac{\eta}{(1+\eps)(1+\eta)}} \left\{ \E (R^{n-1})^{-\frac{1}{\eps}} \right\}^{\frac{\eps}{1+\eps}}
\end{eqnarray}
From the conditions (\ref{hyp1}) or (\ref{hyp2}), we can prove that there exists $\eta > 0$ and $\eps > 0$ s.t.
$$\sup_{n \in \N} \left\{ \E (R^{n-1})^{\frac{1+\eta}{\eta}} \right\}^{\frac{\eta}{(1+\eps)(1+\eta)}} \left\{ \E (R^{n-1})^{-\frac{1}{\eps}} \right\}^{\frac{\eps}{1+\eps}} < + \infty.$$

First assume that (\ref{hyp1}) holds. Then
\begin{eqnarray} \label{210206-5} \nonumber
(R^{n-1})^{-\frac{1}{\eps}} & = & \exp \left[ -\frac{1}{\eps} \int_{0}^{T} g(s,Z^{n-1}_{s}) dW_{s} + \frac{1}{2 \eps} \int_{0}^{T}  g(s,Z^{n-1}_{s})^{2} ds \right] \\ \nonumber
& = & \exp \left[ \int_{0}^{T} \left( - \frac{g(s,Z^{n-1}_{s})}{\eps} \right) dW_{s} - \frac{1}{2} \int_{0}^{T} \left( \frac{g(s,Z^{n-1}_{s})}{\eps} \right)^{2} ds \right] \\
& & \quad \times \exp \left[ \frac{1}{2 \eps^{2}} (1 + \eps) \int_{0}^{T} g(s,Z^{n-1}_{s})^{2} ds  \right]
\end{eqnarray}
Now if
$$\Gamma^{n-1,\eps} = - \int_{0}^{T} \frac{g(u,Z^{n-1}_{u})}{\eps} dW_{u},$$
we have for $C > 1$
$$\E \left[ \exp \left( \frac{\sqrt{C}}{2} \Gamma^{n-1,\eps} \right) \right] \leq \E \left[ \exp \left( \frac{C}{2} \langle \Gamma^{n-1,\eps} \rangle \right) \right]^{1/2} \leq \E \left[ \exp \left( \frac{C}{2 \eps^{2}} \langle L \rangle_{T} \right) \right]^{1/2} < + \infty,$$
when $C/\eps^{2} = \kappa$. Thus
$$\E \left[ \exp \left[ \int_{0}^{T} \left( - \frac{g(s,Z^{n-1}_{s})}{\eps} \right) dW_{s} - \frac{1}{2} \int_{0}^{T} \left( \frac{g(s,Z^{n-1}_{s})}{\eps} \right)^{2} ds \right] \right]^{q} < + \infty$$
when $1/q + 1/p =1$ and
$$\frac{\sqrt{p}}{\sqrt{p} -1 } = C = \eps \sqrt{\kappa} \Longleftrightarrow p = \frac{\kappa \eps^{2}}{\left( \eps \sqrt{\kappa} -1 \right)^{2}}.$$
And we have
$$\E \exp \left[ \frac{p}{2 \eps^{2}} (1 + \eps) \int_{0}^{T} g(s,Z^{n-1}_{s})^{2} ds  \right] \leq \E \exp \left[ \frac{p(1+\eps)}{2 \eps^{2}} \int_{0}^{T} \phi_{s}^{2} ds \right]  < + \infty,$$
if
$$\frac{p(1+\eps)}{\eps^{2}} \leq \kappa \Longleftrightarrow \eps \geq \frac{1+2\sqrt{\kappa}}{\kappa} \Longleftrightarrow 1+ \eps = \frac{\kappa+2 \sqrt{\kappa} +1}{\kappa}.$$
From (\ref{210206-5}) and with H\"{o}lder's inequality we deduce
that $\displaystyle \E R_{n-1}^{-\frac{1}{\eps}} \leq C$. With
(\ref{210206-2}) we already know that there exists $\eta$ s.t.
$\displaystyle \E R_{n-1}^{\frac{1+\eta}{\eta}} \leq C$. We have
to take $\displaystyle \sqrt{1+\eta} =\sqrt{p^*}=
\frac{\sqrt{\kappa}}{\sqrt{\kappa} - 1}.$

Assume that (\ref{hyp2}) holds. Then we already know (\ref{210206-3}): there exists $\eta > 0$ such that
$\displaystyle \E R_{n-1}^{\frac{1+\eta}{\eta}} \leq C$, if $\eta$ satisfies
$$\| L \| < \theta \left( \frac{1+\eta}{\eta} \right).$$
We use theorem 2.4 in \cite{Kaz} in order to prove that
$\displaystyle \E (R^{n-1})^{-\frac{1}{\eps}} \leq C$. We must
choose $\eps$ s.t.
$$\| L \| \leq \sqrt{2} \left( \sqrt{1+\eps} - 1 \right).$$

The two constants $\eta$ and $\eps$ depend on the constant $\kappa$ in (\ref{hyp1}) or the BMO-norm $\| L \|$ in (\ref{hyp2}). Coming back to (\ref{210206-4}) we deduce that:
$$\E \left( \int_{0}^{T} e^{\beta \Phi_s} \left( Z^{n}_{s} \right)^{2} ds \right)^{\frac{p}{2}} \leq C \left\{ \E e^{\beta \frac{p(1+\eps)}{2}(1+\eta) \Phi_T}\left| \xi \right|^{p(1+\eps)(1+\eta)} \right\}^{\frac{1}{(1+\eps)(1+\eta)}}.$$
Remark now that $(1+\eps)(1+\eta) = \Psi$. Thereby from Assumption (\ref{hyp3}), if $\delta > \Psi$, the desired boundedness follows for some $p > 1$ such that $\delta \ge p\Psi$ and by choosing $\beta > 0$ such that $\alpha \geq \beta p \Psi / 2$.
\epf

\begin{pr}
The sequence $(Z_n)_{n \in \N}$ converges in $L^{2}([0,T] \times \prb)$.
\end{pr}
\pf
Applying It\^o's formula we have
\begin{eqnarray*}
&& e^{\beta \Phi_t} |Y^{n+1}_t - Y^{n}_t|^2 + \int_t^T e^{\beta \Phi_u} |Z^{n+1}_u - Z^{n}_u|^2 du = - \beta \int_t^T \phi_u^2 e^{\beta \Phi_u} |Y^{n+1}_u - Y^{n}_u|^2 du \\
&& \qquad -2 \int_t^T e^{\beta \Phi_u} (Y^{n+1}_u-Y^n_u)(-Z^{n+1}_u g(u,Z^{n}_u) + Z^{n}_u g(u,Z^{n-1}_u)) du \\
&& \qquad -2 \int_t^T e^{\beta \Phi_u} (Y^{n+1}_u-Y^n_u)(Z^{n+1}_u - Z^n_u) dW_u \\
&& =  - \beta \int_t^T \phi_u^2 e^{\beta \Phi_u} |Y^{n+1}_u - Y^{n}_u|^2 du + 2\int_t^T e^{\beta \Phi_u} (Y^{n+1}_u-Y^n_u)(Z^{n+1}_u - Z^n_u) g(u,Z^{n}_u) du \\
&& \qquad + 2\int_t^T e^{\beta \Phi_u} (Y^{n+1}_u-Y^n_u)(f(u,Z^{n}_u)-f(u,Z^{n-1}_u)) du \\
&& \qquad + 2\int_t^T e^{\beta \Phi_u} (Y^{n+1}_u-Y^n_u)g(u,Z^{n-1}_u)(Z^{n-1}_u - Z^n_u) du \\
&& \qquad -2 \int_t^T e^{\beta \Phi_u} (Y^{n+1}_u-Y^n_u)(Z^{n+1}_u - Z^n_u) dW_u.
\end{eqnarray*}

Recall that $g$ is bounded by the process $\phi$. Hence with some
positive constants $\eps$ and $\eta$
\begin{eqnarray*}
&&\int_t^T e^{\beta \Phi_u} |Z^{n+1}_u - Z^{n}_u|^2 du \leq \int_t^T \left(\frac{1}{\eps} + 2\frac{1}{\eta} - \beta \right) \phi_u^2 e^{\beta \Phi_u} |Y^{n+1}_u - Y^{n}_u|^2 du \\
&& \qquad + \eps \int_t^T e^{\beta \Phi_u} |Z^{n+1}_u - Z^n_u|^2 du \\
&& \qquad + 2\eta \int_t^T e^{\beta \Phi_u} |Z^{n}_u-Z^{n-1}_u|^2 du \\
&& \qquad -2 \int_t^T e^{\beta \Phi_u} (Y^{n+1}_u-Y^n_u)(Z^{n+1}_u - Z^n_u) dW_u.
\end{eqnarray*}
Choosing $\beta$ such that
\begin{equation} \label{eq1_12_05_08}
\frac{1}{\eps} + 2\frac{1}{\eta} = \beta,
\end{equation}
we have
\begin{eqnarray*}
(1-\eps) \int_0^T e^{\beta \Phi_u} |Z^{n+1}_u - Z^{n}_u|^2 du & \leq & 2 \eta \int_0^T e^{\beta \Phi_u} |Z^{n}_u-Z^{n-1}_u|^2 du \\
& & -2 \int_t^T e^{\beta \Phi_u} (Y^{n+1}_u-Y^n_u)(Z^{n+1}_u - Z^n_u) dW_u.
\end{eqnarray*}
If $\alpha > 4.5 \Psi$, then we can choose $\beta>9$ such that $\alpha \geq \beta \Psi / 2$ (see the end of the proof of Proposition \ref{boundedness}) and $\eps$ and $\eta$ such that (\ref{eq1_12_05_08}) holds with $2 \eta / (1-\eps) < 1$. Since the conclusion of Proposition \ref{boundedness} holds, the local martingale in the previous expression is a true martingale. Hence taking the expectation we obtain:
$$\E \int_0^T e^{\beta \Phi_u} |Z^{n+1}_u - Z^{n}_u|^2 du \leq \frac{2 \eta}{1-\eps} \E \int_0^T e^{\beta \Phi_u} |Z^{n}_u-Z^{n-1}_u|^2 du.$$
Therefore the sequence $(Z_n)_{n \in \N}$ converges in $L^2([0,T] \times \prb)$.
\epf

\begin{lem} \label{aeconv}
There exists a subsequence of $Z^n$ (still denoted $Z^n$) which converges $\prb \otimes \lambda$-a.e. to some process $Z$.
\end{lem}

\begin{lem}
The sequence $R^n_T$ converges also $\prb$-a.s. to
$$R_T = \exp \left( \int_0^T g(s,Z_s) dW_s - \frac{1}{2} \int_0^T (g(s,Z_s))^2 ds \right).$$
\end{lem}
\pf
We may w.l.o.g. assume that $g(s,.)$ is continuous for all $s \in [0,T]$. The rest follows from Lemma \ref{aeconv}.
\epf

Equipped with these results, we are now in a position to state our existence Theorem.

\begin{thm}\label{existence_general}
Suppose Assumption (H1) holds. There exists a probability measure
$\qprb$ equivalent to $\prb$ and an adapted process $Z$ such that
$\E \int_0^T |Z_s|^{2} ds <\infty$ such that, setting
$$R_T = \exp \left( \int_0^T g(s, Z_s) d W_s - \frac{1}{2} \int_0^T g(s, Z_s)^2 ds \right),\quad W^{\qprb} = W - \int_0^\cdot g(s, Z_s) ds,$$
we have
$$ \qprb = R_T \cdot \prb,$$
and such that the pair $(Y, Z)$ defined by
$$Y = \E^{\qprb}(\xi|\tri_\cdot) = \E^{\qprb}(\xi) + \int_0^\cdot Z_s d W^{\qprb}_s$$
solves the BSDE (\ref{BSDE}).
\end{thm}

\pf Using Theorem \ref{tightness_lipschitz}, choose a probability
measure $\qprb$ and another subsequence of the corresponding
subsequence of $(\mathbb{Q}_n)_{n\in\N}$ which converges weakly to
$\qprb$. We denote this subsequence again by
$(\mathbb{Q}_n)_{n\in\N}$ and the corresponding subsequence of
controls by $(Z^n)_{n\in\N}$. We have:
$$\qprb = R_T \cdot \prb.$$
Moreover for all $n \in \N$,
$$Y^n_t = \E^{n-1}(\xi) + \int_0^t Z^n_s d W^n_s = \E^{n-1}(\xi) + \int_0^t Z^n_s d W_s - \int_0^t Z^n_s g(s,Z^{n-1}_s) ds.$$
The only thing we have to prove, is that the sequence $Y^n_0 = \E^{n}(\xi)$ also converges. But $Y^n_0 = \E^{n}(\xi) = \E (\xi R^n)$, and $\xi$ belongs to $L^{\gamma}$, $R^n$ also belongs to some $L^p$ space with $1/p + 1/\gamma = 1$ if and only if
$$\gamma \ge \frac{\kappa}{(\sqrt{\kappa}-1)^2}.$$
But it is true since $\gamma \ge \Psi(\kappa)$. Taking a
subsequence if necessary, we deduce that $Y^n_0$ converges to
$\E^{\qprb} (\xi)$.

Hence we obtain
$$Y_t = \E^{\qprb}(\xi | \tri_\cdot) = \E^{\qprb}(\xi) + \int_0^t Z_s d W^{\qprb}_s,$$
where $W^{\qprb}$ is a $\qprb$-Brownian motion. Finally $(Y, Z)$ solves the BSDE (\ref{BSDE}).
\epf


\end{document}